\documentclass{amsart}

\newtheorem{theorem}{Theorem}[section]
\newtheorem{lemma}[theorem]{Lemma}

\newtheorem{corollary}[theorem]{Corollary}
\newtheorem{proposition}[theorem]{Proposition}
\usepackage{graphicx}

\theoremstyle{definition}

\theoremstyle{remark}

\numberwithin{equation}{section}

\begin{document}

\title[Distance between toroidal surgeries]{Distance between toroidal surgeries on hyperbolic knots in the $3$-sphere}

\author{Masakazu Teragaito}
\address{Department of Mathematics and Mathematics Education, Hiroshima University,
1-1-1 Kagamiyama, Higashi-hiroshima, Japan 739-8524}
\email{teragai@hiroshima-u.ac.jp}
\thanks{
Partially supported by Japan Society for the Promotion of Science,
Grant-in-Aid for Scientific Research (C), 14540082.
}%





\begin{abstract}
For a hyperbolic knot in the $3$-sphere, at most finitely many Dehn surgeries
yield non-hyperbolic $3$-manifolds.
As a typical case of such an exceptional surgery, a toroidal surgery
is one that yields a closed $3$-manifold containing an incompressible torus.
The slope corresponding to a toroidal surgery, called a toroidal slope,
is known to be integral or half-integral.
We show that the distance between two integral toroidal slopes for a hyperbolic knot, except the figure-eight knot,
is at most four.
\end{abstract}

\maketitle

\section{Introduction}

Let $K$ be a knot in the $3$-sphere $S^3$ and let $E(K)=S^3-\mathrm{Int}\,N(K)$ be its exterior.
A \textit{slope\/} is the isotopy class of an essential simple closed curve on $\partial E(K)$.
Then the set of slopes is parameterized by $\mathbf{Q}\cup \{1/0\}$ so that $1/0$ is the meridian slope as in the usual way
(see \cite{R}).
For two slopes $\alpha$ and $\beta$, the distance $\Delta(\alpha,\beta)$ between $\alpha$ and $\beta$ is defined to be
their minimal geometric intersection number.
A slope $m/n$ is called \textit{integral\/} if $|n|=1$, and \textit{half-integral\/} if $|n|=2$.
In other words, an integral slope runs once along the knot, and a half-integral slope runs twice along the knot.

We denote by $K(\alpha)$ the closed $3$-manifold obtained by \textit{$\alpha$-Dehn surgery\/} on $K$, that is,
attaching a solid torus $V$ to $E(K)$ along $\partial E(K)$ in such a way that the slope $\alpha$ bounds a disk in $V$.
A surgery (or slope) is said to be \textit{toroidal\/} if the resulting manifold contains an incompressible torus.
Thurston showed that if $K$ is a hyperbolic knot,
then $K(\alpha)$ is hyperbolic for all but finitely many slopes $\alpha$ \cite{Th}.
If $K(\alpha)$ is not hyperbolic, then it is either reducible, or an atoroidal Seifert fibered manifold, or toroidal,
or a counterexample to the Geometrization Conjecture \cite{Th}.
We focus on the third case.
It is known that if $\alpha$ is a toroidal slope for a hyperbolic knot, then
$\alpha$ is integral or half-integral \cite{GL,GL2}.
There are many examples of integral toroidal surgery, and Eudave-Mu\~{n}oz \cite{EM}
constructed an infinite family of hyperbolic knots $k(\ell,m,n,p)$ admitting
half-integral toroidal surgeries. Recently, Gordon and Luecke \cite{GL3} proved that
the Eudave-Mu\~{n}oz knots $k(\ell,m,n,p)$ are the only hyperbolic knots with half-integral toroidal surgeries.

In this paper, we consider the distance between toroidal slopes on a hyperbolic knot in $S^3$.
The figure-eight knot admits exactly three toroidal slopes $0, 4$ and $-4$ \cite{Th}.
Note that $\Delta(-4,4)=8$.
If a hyperbolic knot is not the figure-eight knot, then
the distance between two toroidal slopes is at most $5$ by Gordon \cite{G}.
(There are exactly four hyperbolic $3$-manifolds which admit two toroidal slopes with distance at least $6$. 
They all are obtained from the Whitehead link by some Dehn surgery on one component.
Among those, only the figure-eight knot exterior can be embedded in $S^3$ by homological reason.)
This upper bound $5$ is sharp.
For example, the Eudave-Mu\~{n}oz knot $k(2,-1,n,0)$ ($n\ne 1$) admits two toroidal slopes $25n-\frac{37}{2}$ and $25n-16$ as shown in \cite{EM2}, where
$\Delta(25n-\frac{37}{2},25n-16)=5$. (When $n=1$, $k(2,-1,1,0)$ is the trefoil.)
Notice that one slope is half-integral, and that
$k(2,-1,0,0)$ is the $(2,-3,-7)$-pretzel knot.
The purpose of this paper is to show that
we can reduce the upper bound when both of toroidal slopes are integral.

\begin{theorem}\label{main1}
Let $K$ be a hyperbolic knot in $S^3$, which is not the figure-eight knot.
If $\alpha$ and $\beta$ are two integral toroidal slopes for $K$, then
$\Delta(\alpha,\beta)\le 4$.
\end{theorem}

This is sharp.
For example, the twist knot $C[2n,2]$ in Conway's notation with $n\ge 1$ admits two integral toroidal slopes $0$ and $4$ \cite{BW}.
Although it may be too optimistic, there is a possibility that only twist knots and the $(-2,3,7)$-pretzel knot
admit two integral toroidal slopes with distance  $4$.
%
%
%
%

\begin{corollary}\label{main2}
If a hyperbolic knot $K$ in $S^3$ admits two toroidal slopes $\alpha$ and $\beta$ with $\Delta(\alpha,\beta)=5$, then
$K$ is an Eudave-Mu\~{n}oz knot.
\end{corollary}

\begin{proof}
By Theorem \ref{main1}, one of $\alpha, \beta$ is half-integral.
Then $K$ is an Eudave-Mu\~{n}oz knot by \cite{GL3}.
\end{proof}

Among Eudave-Mu\~{n}oz knots, only $k(2,-1,n,0)$ $(n\ne 1)$ seems to admit two toroidal slopes with distance $5$.
But this is still an open question.


It is conjectured that a hyperbolic knot in $S^3$ admits at most three toroidal slopes \cite{EM} (see also \cite[Problem 1.77 A(5)]{K}).
%
%
%
%
Our main theorem also gives an upper bound for the number of toroidal slopes.

\begin{corollary}
A hyperbolic knot in $S^3$ admits at most $5$ toroidal slopes.
\end{corollary}

\begin{proof}
Any Eudave-Mu\~{n}oz knot $k=k(\ell,m,n,p)$ admits at least three non-trivial exceptional slopes $s-1,s-\frac{1}{2},s$, where $s$ is an integer
determined by $k$ \cite{EM2}.
In fact, $s-\frac{1}{2}$ is the only half-integral toroidal slope for $k$ \cite{GWZ}, and $s-1$ and $s$ yield atoroidal Seifert fibered manifolds.
Since the distance between two toroidal slopes is at most $5$ by \cite{G},
the only possible toroidal slopes for $k$ are $s-3,s-2,s-\frac{1}{2},s+1,s+2$.
But both of $s-3$ and $s+2$ cannot be toroidal by Theorem \ref{main1}.
Thus $k$ admits at most $4$ toroidal slopes.

The figure-eight knot admits three toroidal slopes as stated before.
For the other hyperbolic knots, any toroidal slope is integral and the distance between such two slopes is at most $4$ by Theorem \ref{main1}.
Hence there are at most $5$ toroidal slopes.
\end{proof}

In Section \ref{preliminaries}, we prepare the basic tool, a pair of labelled graphs, to show Theorem \ref{main1}.
Also, some fundamental properties are shown there.
Section \ref{reduced} is devoted to examine a reduced graph supported in a disk or an annulus.
The results will be used in Section \ref{generic}.
Sections \ref{s2t4} and \ref{s2t2} will treat two special cases, and the proof of Theorem \ref{main1} will be completed.
In the last section, we propose some questions concerning toroidal slopes of hyperbolic knots.

\section{Preliminaries}\label{preliminaries}

Throughout this paper, we fix a hyperbolic knot $K$ in $S^3$, which is not the figure-eight knot.
Let $\alpha$ and $\beta$ be two integral toroidal slopes for $K$.
By \cite{G}, $\Delta(\alpha,\beta)\le 5$.
We assume $\Delta(\alpha,\beta)=5$ to prove Theorem \ref{main1}.
Note that $K(\alpha)$ and $K(\beta)$ are irreducible under this assumption.
(For, if $K(\alpha)$ is reducible, then $\Delta(\alpha,\beta)\le 3$ by \cite{O, Wu}. Similarly for $K(\beta)$.)
Let $\widehat{S}$ be an incompressible torus in $K(\alpha)$.
We may assume that $\widehat{S}$ meets the attached solid torus $V_\alpha$ in $s$ meridian disks
$u_1,u_2,\dots,u_s$, numbered successively along $V_\alpha$, and
that $s$ is minimal over all choices of $\widehat{S}$.
Let $S=\widehat{S}\cap E(K)$.
Then $S$ is a punctured torus properly embedded in $E(K)$ with $s$ boundary components $\partial_iS=\partial u_i$,
each of which has slope $\alpha$.
By the minimality of $s$, $S$ is incompressible, and then boundary-incompressible in $E(K)$.
Similarly, we choose an incompressible torus $\widehat{T}$ in $K(\beta)$ which intersects the
attached solid torus $V_\beta$ in $t$ meridian disks $v_1,v_2,\dots,v_t$, numbered successively along $V_\beta$, where $t$ is minimal as above.
Then we have another incompressible and boundary-incompressible punctured torus $T=\widehat{T}\cap E(K)$,
which has $t$ boundary components $\partial_jT=\partial v_j$.

\begin{proposition}
The genus of $K$ is not one.
\end{proposition}

\begin{proof}
Assume that $K$ has genus one.
By \cite{Te2}, if $r$ is an integral toroidal slope for $K$, then $|r|=0,1,2$ or $4$.
Furthermore, if $|r|=4$ then $K$ is a twist knot.
Since $\Delta(\alpha,\beta)=5$, either slope is $-4$ or $4$.
Thus $K$ is a twist knot.
But a twist knot does not admit two toroidal slopes with distance $5$ \cite{BW}.
\end{proof}

\begin{lemma}
$\widehat{S}$ and $\widehat{T}$ are separating.
\end{lemma}

\begin{proof}
Assume $\widehat{S}$ is non-separating.
Then $\alpha=0$ by homological reason.
Thus $K(0)$ contains a non-separating torus $\widehat{S}$.
But this implies that $K$ has genus one \cite[Corollary 8.3]{Ga}.
Similarly for $\widehat{T}$.
\end{proof}

Thus $s$ and $t$ are non-zero and even.

We may assume that $S$ intersects $T$ transversely.
Then $S\cap T$ consists of arcs and circles.
Since both surfaces are incompressible, we can assume that no circle component of $S\cap T$
bounds a disk in $S$ or $T$.
Furthermore, it can be assumed that $\partial_iS$ meets $\partial_jT$ in $5$ points for any pair of $i$ and $j$.

\begin{lemma}\label{jumping}
Let $a_1,a_2,a_3,a_4,a_5$ be the points of $\partial_iS\cap\partial_jT$, numbered so that they appear successively on $\partial_iS$.
Then $a_1,a_2,a_3,a_4,a_5$ also appear successively on $\partial_jT$.
In particular, two points of $\partial_iS\cap\partial_jT$ are successive on $\partial_iS$ if and only if they are successive in $\partial_jT$.
\end{lemma}

\begin{proof}
This immediately follows from that both slopes $\alpha$ and $\beta$ are integral.
\end{proof}

In the literature, for example \cite{GW}, this fact is stated that the jumping number is one.

Let $G_S$ be the graph on $\widehat{S}$ consisting of the $u_i$ as (fat) vertices, and the arc components
of $S\cap T$ as edges.
Define $G_T$ on $\widehat{T}$ similarly.
Throughout the paper, two graphs on a surface are considered to be equivalent if there is a homeomorphism of the surface
carrying one graph to the other.
Note that $G_S$ and $G_T$ have no trivial loops.

For an edge $e$ of $G_S$ incident to $u_i$,
the endpoint of $e$ is labelled $j$ if it is in $\partial u_i\cap \partial v_j=\partial_iS\cap\partial_jT$.
Similarly, label the endpoints of each edge of $G_T$.
Thus the labels $1,2,\dots,t$ ($1,2,\dots,s$, resp.) appear in order around each vertex of $G_S$ ($G_T$, resp.) repeated $5$ times.
Each vertex $u_i$ of $G_S$ has degree $5t$, and each $v_j$ of $G_T$ has degree $5s$.

Let $G=G_S$ or $G_T$.
Two vertices of $G$ are said to be \textit{parallel\/} if their numbers have the same parity, 
otherwise they are \textit{antiparallel\/}.
An edge of $G$ is a \textit{positive\/} edge if it connects parallel vertices.
Otherwise it is a \textit{negative\/} edge.
Possibly, a positive edge is a loop.
An endpoint of a positive (negative, resp.) edge around a vertex is called a \textit{positive} (\textit{negative\/}, resp.) \textit{edge endpoint\/}.

\begin{lemma}[The parity rule]
An edge $e$ is positive in a graph if and only if it is negative in the other graph.
\end{lemma}

\begin{proof}
This can be found in \cite{CGLS}.
\end{proof}

\begin{lemma}\label{bothparallel}
There is no pair of edges which are parallel in both graphs.
\end{lemma}

\begin{proof}
This is \cite[Lemma 2.1]{G}.
\end{proof}

If an edge $e$ of $G_S$ is incident to $u_i$ with label $j$, the it is called a \textit{$j$-edge at} $u_i$.
Then $e$ is also an $i$-edge at $v_j$ in $G_T$.
If $e$ has labels $j_1,j_2$ at its endpoints, then $e$ is called a \textit{$\{j_1,j_2\}$-edge}, or $j_1j_2$-edge.

A cycle in $G$ consisting of positive edges is a \textit{Scharlemann cycle\/} if it bounds a disk face of $G$ and all edges
in the cycle are $\{i,i+1\}$-edges for some label $i$.
The number of edges in a Scharlemann cycle is called the \textit{length\/} of the Scharlemann cycle, and
the set $\{i,i+1\}$ is called its \textit{label pair}.
A Scharlemann cycle of length two is called an \textit{$S$-cycle} in short.

Let $e_1,e_2,e_3,e_4$ be the four parallel positive edges in $G$ numbered in order.
If $G$ has at least four labels, and the middle two edges $e_2$ and $e_3$ form an $S$-cycle,
then the cycle defined by $e_1$ and $e_4$ is called an \textit{extended $S$-cycle}.
(There is a notion of an extended Scharlemann cycle of arbitrary length, but this is enough for our purpose.)

\begin{lemma}\label{Schsupp}
Let $\rho$ be a Scharlemann cycle in $G_S$.
Then the edges of $\rho$ cannot lie in a disk in $\widehat{T}$.
This also holds for Scharlemann cycles in $G_T$.
\end{lemma}

\begin{proof}
Without loss of generality, assume that $\rho$ has label pair $\{1,2\}$.
Let $f$ be the face of $G_S$ bounded by $\rho$, and let $V_{12}$ be the part of the attached solid torus
$V_\beta$ between two meridian disks $v_1$ and $v_2$.
(When $t=2$, choose one such that $\partial f$ runs on $\partial V_{12}-v_1\cup v_2$.)
Assume that the edges of $\rho$ lie in a disk $D$ in $\widehat{T}$.
Then $N(D\cup V_{12}\cup f)$ is a lens space minus an open $3$-ball.
Since $K(\beta)$ is irreducible, $K(\beta)$ is a lens space.
But this contradicts that $K(\beta)$ is toroidal.
\end{proof}

For the rest of this section, let $G=G_S$ or $G_T$.
Assume that $G$ has $q$ labels.

\begin{lemma}\label{Scharlemann}
Let $q\ge 4$.
\begin{itemize}
\item[(1)] $G$ cannot contain an extended $S$-cycle.
\item[(2)] If $q=4$, then $G$ cannot contain two $S$-cycles with disjoint label pairs.
\item[(3)] $G$ cannot contain three $S$-cycles with mutually disjoint label pairs.
\item[(4)] If there are two $S$-cycles with disjoint label pairs $\{i,i+1\}$ and $\{j,j+1\}$, then $i$ and $j$ have the same parity.
\end{itemize}
\end{lemma}

\begin{proof}
(1) is \cite[Lemma 2.10]{BZ2}.
(2) For convenience, we assume $G=G_S$.
We can assume that $\rho_1$ and $\rho_2$ are $S$-cycles with label pairs $\{1,2\}$ and $\{3,4\}$, respectively.
Let $f_i$ be the face of $G_S$ bounded by $\rho_i$, $i=1,2$.
Denote by $V_{12}$ and $V_{34}$ the parts of the attached solid torus $V_\beta$ lying between two meridian disks $v_1$ and $v_2$, and $v_3$ and $v_4$, 
respectively.
Then shrinking $V_{12}$ radially to its core in $V_{12}\cup f_1$ gives a M\"{o}bius band $B_1$ such that $\partial B_1$ is the loop
on $\widehat{T}$ formed by the edges of $\rho_1$.
Similarly, $V_{34}\cup f_2$ gives another M\"{o}bius band $B_2$ whose boundary is disjoint from $\partial B_1$.
Let $A$ be an annulus between $\partial B_1$ and $\partial B_2$ on $\widehat{T}$.
Then $B_1\cup A\cup B_2$ is a Klein bottle $\widehat{F}$ in $K(\beta)$, which meets $V_\beta$ in two meridian disks (after a perturbation).
Then $F=\widehat{F}\cap E(K)$ gives a twice-punctured Klein bottle in $E(K)$.
By attaching a suitable annulus on $\partial E(K)$ to $F$ along their boundaries, 
we have a closed non-orientable surface in $S^3$, a contradiction.
(3) and (4) are \cite[Lemma 1.10]{Wu} and \cite[Lemma 1.7]{Wu}, respectively.
\end{proof}

Let $e_1,e_2,\dots,e_q$ be $q$ mutually parallel negative edges in $G$ numbered successively, each connecting vertex $x$ to $y$.
Suppose that $e_i$ has label $i$ at $x$ for $1\le i\le q$.
Then this family defines a permutation $\sigma$ of the set $\{1,2,\dots,q\}$ such that
$e_i$ has label $\sigma(i)$ at $y$.
In fact, $\sigma(i)\equiv i+r \pmod{q}$ for some even $r$.
We call $\sigma$ the \textit{permutation associated to the family\/}. 
It is well-defined up to inversion.
Note that $\sigma$ has at least two orbits by the parity rule, and all orbits of $\sigma$ have the same length.

\begin{lemma}\label{G}
Let $q\ge 4$.
\begin{itemize}
\item[(1)] Any family of parallel positive edges in $G$ contains at most $q/2+1$ edges.
Moreover, if the family contains $q/2+1$ edges, then two adjacent edges on one end form an $S$-cycle.
\item[(2)] Any family of parallel negative edges in $G$ contains at most $q$ edges.
\end{itemize}
\end{lemma}

\begin{figure}[tb]
\includegraphics*[scale=0.4]{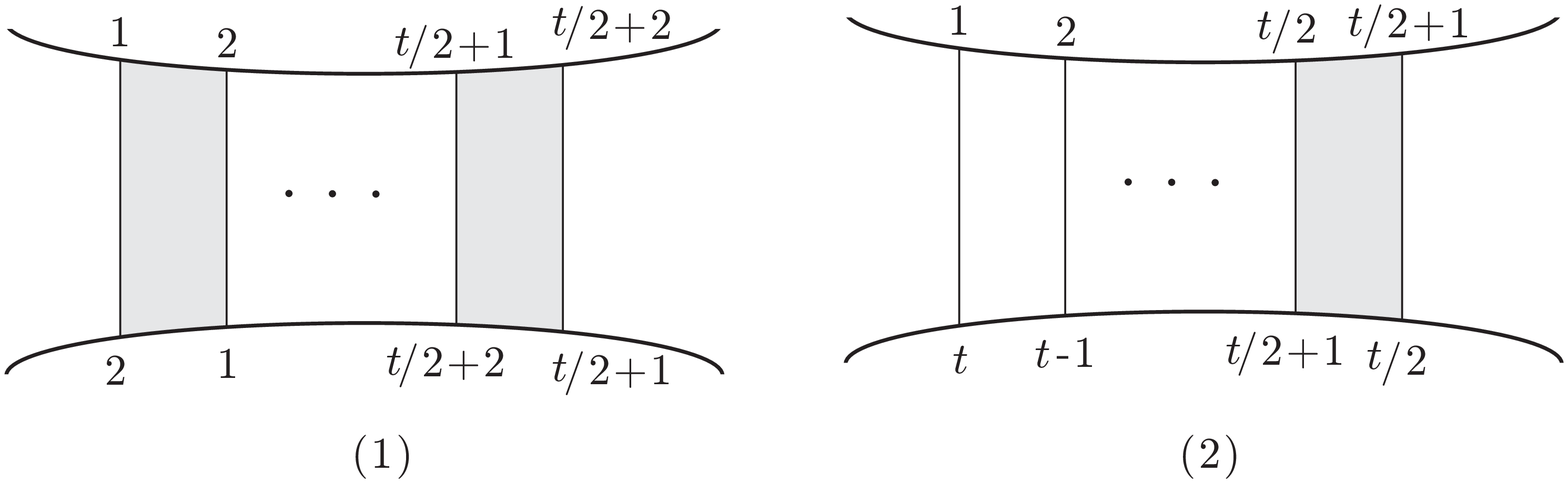}
\caption{}\label{positive}
\end{figure}

\begin{proof}
For convenience, let $G=G_S$.  

(1) Such a family contains at most $t/2+2$ edges \cite[Lemma 1.4]{Wu}.
Furthermore, if it contains $t/2+2$ edges, then the labels can be assumed as in Figure \ref{positive}(1).
Then there are two $S$-cycles with disjoint label pairs $\{1,2\}$ and $\{t/2+1,t/2+2\}$.
By the same construction as in the proof of Lemma \ref{Scharlemann}(2), we have two M\"{o}bius bands $B_1$ and $B_2$ and an annulus $A$ on $\widehat{T}$
as before.
In $\widehat{T}$, two vertices $v_{i}$ and $v_{t-i+3}$ are connected by an edge in the family for $i=3,4,\dots,t/2$.
Hence $\mathrm{int}\,A$ contains an even number of vertices.
Then $B_1\cup A\cup B_2$ is a Klein bottle which meets $V_\beta$ in an even number of meridian disks.
This leads to a contradiction as before.

If the family contains $t/2+1$ edges, then it contains an $S$-cycle \cite[Lemma 2.6.6]{CGLS}.
Since there is no extended $S$-cycle by Lemma \ref{Scharlemann}(1), the last two edges form an $S$-cycle as shown in Figure \ref{positive}(2),
up to relabeling.

(2) Let $\{e_1,e_2,\dots,e_t,e_1'\}$ be a family of $t+1$ parallel negative edges, connecting 
vertex $u_j$ to $u_k$.
We can assume that $e_i$ has label $i$ at $u_j$ for $1\le i\le t$ and $e_1'$ has label $1$ at $u_j$.
Let $\sigma$ be the permutation associated to the family $\{e_1,e_2,\dots,e_t\}$.
Thus $e_i$ has label $\sigma(i)$ at $u_k$.
For an orbit $\theta$ of $\sigma$, let $C_\theta$ be the cycle in $G_T$ consisting of $e_i$ for $i\in \theta$.
Then $C_\theta$ does not bound a disk on $\widehat{T}$ \cite[Lemma 2.3]{G}.
Since $i$ and $\sigma(i)$ have the same parity by the parity rule,
$\sigma$ has at least two orbits.
Thus all cycles corresponding to the orbits of $\sigma$ are essential and mutually parallel on $\widehat{T}$.
Let $C_1$ be the cycle corresponding to the orbit containing $1$.

Now, consider $e_1'$.
Then $e_1'$ connects $v_1$ and $v_{\sigma(1)}$ in $G_T$.
By Lemma \ref{bothparallel}, it is not parallel to $e_1$ in $G_T$.
Because of the existence of another cycle, the new cycle obtained from $C_1$ by exchanging $e_1$ by $e_1'$ bounds
a disk in $\widehat{T}$.
But this is impossible by \cite[Lemma 2.3]{G} again.
\end{proof}

\begin{lemma}\label{consec}
If $q=4$, then there are no consecutive $4$ families of $q/2+1$ parallel positive edges at any vertex of $G$.
\end{lemma}

\begin{proof}
If there are such $4$ families, then there are two $S$-cycles with disjoint label pairs among those, which is
impossible by Lemma \ref{Scharlemann}(2).
\end{proof}

\section{Reduced graphs on tori}\label{reduced}

Let $G=G_S$ or $G_T$, and let $F$ be the surface where $G$ lies.
The \textit{reduced graph\/} $\overline{G}$ of $G$ is the graph obtained from $G$ by amalgamating each family of mutually parallel edges into
a single edge.
Let $G^+$ be the subgraph of $G$ consisting of all vertices  and all positive edges of $G$. 
Then $\overline{G}^+$ is also defined similarly.
In this section, we examine the reduced graphs $\overline{G}$ and $\overline{G}^+$.
The results here will be used in the next section.

Let $\Gamma$ be a component of $\overline{G}^+$.
If there is a disk $D$ in $F$ such that $\mathrm{int}\,D$ contains $\Gamma$, then we say that $\Gamma$ has a \textit{disk support}.
Also, if there is an annulus $A$ in $F$ such that $\mathrm{int}\,A$ contains $\Gamma$ and $\Gamma$ does not have a disk support,
then we say that $\Gamma$ has an \textit{annulus support}.
If $\Gamma$ has neither a disk nor an annulus support, then we say that $\Gamma$ has a \textit{torus support}.

Now, suppose that $\Gamma$ has a support $E$, where $E$ is a disk or an annulus.
A vertex $x$ of $\Gamma$ is called an \textit{outer vertex\/} if there is an arc $\xi$ connecting $x$ to $\partial E$
whose interior is disjoint from $\Gamma$.
Define an \textit{outer edge\/} similarly.
Then $\partial \Gamma$ denotes the subgraph of $\Gamma$ consisting of all outer vertices and all outer edges of $\Gamma$.
A vertex $x$ of $\Gamma$ is called a \textit{cut vertex\/} if $\Gamma-x$ has more components than $\Gamma$.

Suppose that $\Gamma$ has an annulus support $A$.
A vertex $x$ of $\Gamma$ is a \textit{pinched vertex\/} if
there is a spanning arc of $A$ which meets $\Gamma$ in only $x$.
An edge $e$ of $\Gamma$ is a \textit{pinched edge\/} if
there is a spanning arc of $A$ which meets $\Gamma$ in only one point on $e$.
Clearly, both endpoints of a pinched edge are pinched vertices.

We say that $\Gamma$ is an \textit{extremal component\/} of $\overline{G}^+$ if $\Gamma$ has a support which is disjoint from
the other components of $\overline{G}^+$.
Remark that $\overline{G}^+$ has at least two components, because $G$ has vertices of distinct parities.

\begin{lemma}\label{DorA}
$\overline{G}^+$ has an extremal component with a disk support or an annulus support.
\end{lemma}

\begin{proof}
There are only three possibilities for the support of a component; a disk, an annulus, or a torus.
If there is a component with a torus support, then any other component has an disk support. 
The conclusion immediately follows from this.
\end{proof}

Let $x$ be a vertex of $G$.
Then $x$ is called an \textit{interior vertex\/} if there is no negative edge incident to $x$ in $G$.
Since $\overline{G}$ and $\overline{G}^+$  have the same vertex set as $G$,
we may call a vertex of $\overline{G}$ or $\overline{G}^+$ an interior vertex when it is an interior vertex of $G$.
In particular, if $x$ is in an extremal component of $\overline{G}^+$ with a disk or an annulus support, and
it is not an outer vertex, then $x$ is an interior vertex.

A vertex $x$ is said to be \textit{good\/} if all positive edge endpoints around $x$ are successive in $G$.
Thus an interior vertex is good.
When $x$ is a vertex of an extremal component $\Gamma$ of $\overline{G}^+$ with a disk or an annulus support,
$x$ is good if
\begin{itemize}
\item[(i)] $x$ is not a cut vertex of $\Gamma$ if $\Gamma$ has a disk support; or
\item[(ii)] $x$ is neither a cut vertex nor a pinched vertex of $\Gamma$ if $\Gamma$ has an annulus support.
\end{itemize}

\begin{lemma}\label{Dsupp}
Let $\Gamma$ be an extremal component of $\overline{G}^+$.
Assume that $\Gamma$ has a disk support, and that $\Gamma$ is not a single vertex.
\begin{itemize}
\item[(1)] If each interior vertex of $\Gamma$ has degree at least $6$, then $\Gamma$ has two good vertices of degree at most $3$.
\item[(2)] If $\Gamma$ has no interior vertex, then $\Gamma$ has two good vertices of degree at most $2$.
\end{itemize}
\end{lemma}

\begin{proof}
These are \cite[Lemma 2.3]{Wu} and \cite[Lemma 3.2]{Wu}.
\end{proof}

\begin{lemma}\label{Asupp}
Let $\Gamma$ be an extremal component of $\overline{G}^+$.
Assume that $\Gamma$ has an annulus support, and that $\Gamma$ is not a cycle.
If each interior vertex of $\Gamma$ has degree at least $6$, then $\Gamma$ has a good vertex of degree at most $4$.
\end{lemma}

\begin{proof}
First, consider the case that $\Gamma$ has no cut vertices.

Assume that $\Gamma$ has no pinched vertex.
If any vertex on $\partial \Gamma$ has degree at least $5$, then
take a double of $\Gamma$ along two boundary cycles.
Then we have a graph in a torus, whose vertices have degree at least $6$, and
at least two vertices, coming from the vertices on $\partial\Gamma$,
have degree at least $8$.
Also, the graph has no trivial loop and parallel edges.
This is impossible by an Euler characteristic argument.
Therefore, some vertex on $\partial \Gamma$ has degree at most $4$.
Clearly, it is a good vertex.

Assume that $\Gamma$ has a pinched vertex.
Contract each pinched edge into a point if necessary.
Let $\Gamma'$ be the resulting graph.
By our assumption that $\Lambda$ is not a cycle, neither is $\Gamma'$ .
Moreover, any pinched vertex of $\Gamma'$ has degree at least $4$.

If there is only one pinched vertex $x$,
then split $\Gamma'$ at $x$ to obtain $\Gamma''$ having a disk support.
Then $x$ splits into $x_1$ and $x_2$.
Let $V, E, F$ be the number of vertices, edges and faces of $\Gamma''$ as a graph
in a disk.
Let $X$ be the number of vertices on the boundary of $\Gamma''$ except $x_1$ and $x_2$.
Assume that those $X$ vertices on $\partial \Gamma''$
have degree at least $5$ for a contradiction.
Then  $1=V-E+F$, $5X+ \mathrm{deg}(x_1)+\mathrm{deg}(x_2)+6(V-X-2)\le 2E$, and $3F+(X+2)\le 2E$.
(We use $\mathrm{deg(-)}$ to denote the degree of a vertex.)
Thus we have $\mathrm{deg}(x_1)+\mathrm{deg}(x_2)\le 2$.
Since $\mathrm{deg}(x_1)+\mathrm{deg}(x_2)\ge 4$, this is a contradiction.
Hence we see that $\Gamma$ has a vertex of degree at most $4$ on $\partial \Gamma$,
which is not a pinched vertex.

If there are more than one pinched vertices in $\Gamma'$,
then consider two consecutive pinched vertices $y$ and $z$.
Let $\Lambda$ be a subgraph of $\Gamma'$ between $y$ and $z$, which
contains no other pinched vertex.
Then the same argument as above gives 
a desired vertex on $\partial\Gamma$.

Next,
assume that $\Gamma$ has a cut vertex $x$.
Let $\Gamma_1,\Gamma_2,\dots,\Gamma_k$ be the components after splitting along $x$.
If some $\Gamma_i$ has a disk support, then 
$\Gamma_i$ has two good vertices of degree at most $3$ by Lemma \ref{Dsupp},
and one of which is not $x$.
Thus we have a desired vertex.
Otherwise, each $\Gamma_i$ has an annulus support.
Then either $\Gamma$ has no pinched vertex, or $x$ is the unique pinched vertex.
In the former, 
taking a double of $\Gamma$ along two boundary cycles gives
a contradiction as before, unless some vertex on $\partial \Gamma$ has degree at most $4$.
In the latter, split $\Gamma$ at $x$ along a spanning arc of the annulus support meeting $\Gamma$ in only $x$.
Then the same calculation as above gives the conclusion.
\end{proof}

Of course, the conclusion of this lemma is true when $\Gamma$ has no interior vertex.

\begin{proposition}\label{DandAsupp}
If each interior vertex of $\overline{G}^+$ has degree at least $6$, then $\overline{G}^+$ has a vertex of degree at most $4$.
\end{proposition}

\begin{proof}
By Lemma \ref{DorA}, $\overline{G}^+$ has an extremal component $\Gamma$ with a disk or an annulus support.
If $\Gamma$ is a single vertex or a cycle, then the result is obvious.
Otherwise, it follows from Lemmas \ref{Dsupp} and \ref{Asupp}.
\end{proof}

\begin{figure}[tb]
\includegraphics*[scale=0.26]{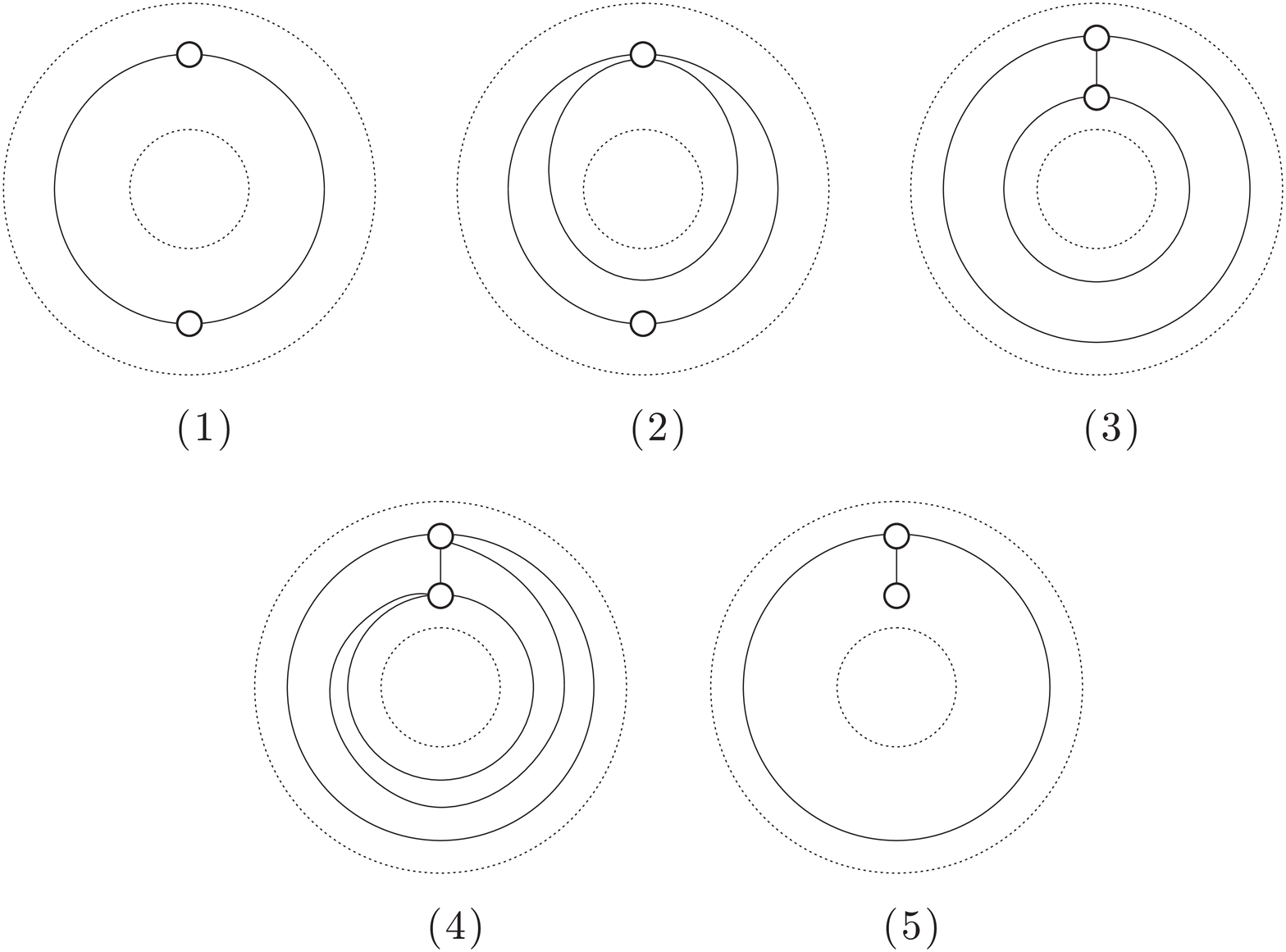}
\caption{}\label{vertex2}
\end{figure}

\begin{lemma}\label{2vertex}
Let $\Gamma$ be a component of $\overline{G}^+$ with an annulus support.
If $\Gamma$ has just two vertices and no interior vertex, then there are five possibilities for $\Gamma$ as shown in Figure \ref{vertex2}. 
\end{lemma}

\begin{proof}
If both vertices are incident to loops, then we have (3) or (4).
If only one vertex is incident to a loop, then $\Gamma$ is (2) or (5).
Finally, if there is no loop, then $\Gamma$ is (1).
\end{proof}

\begin{lemma}\label{3vertex}
Let $\Gamma$ be a component of $\overline{G}^+$ with an annulus support.
Assume that $\Gamma$ is not a cycle.
If $\Gamma$ has just three vertices and no interior vertex,
then $\Gamma$ has a good vertex of degree at most $3$.
\end{lemma}

\begin{proof}
If $\Gamma$ has a block with a disk support, then Lemma \ref{Dsupp} can be applied to the block and
we have a good vertex of degree at most two.
Otherwise, $\partial \Gamma$ consists of two disjoint cycles, or $\Gamma$ has a pinched vertex.
Also, $\Gamma$ has no vertex of degree one.

\begin{figure}[tb]
\includegraphics*[scale=0.28]{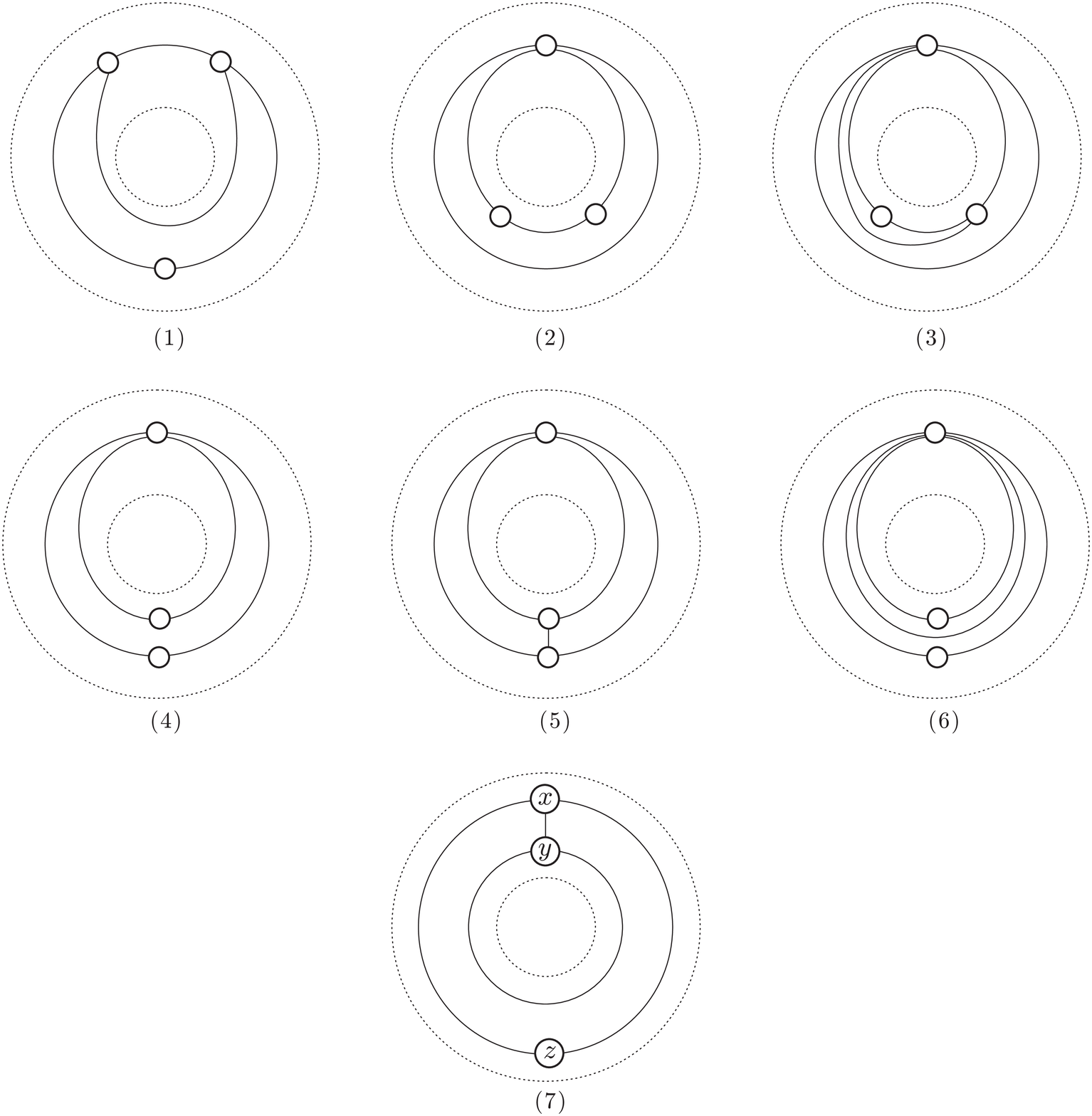}
\caption{}\label{vertex3}
\end{figure}

If $\Gamma$ has two pinched vertices,
there is a pinched edge.
Then $\Gamma$ is either a cycle or the graph as shown in Figure \ref{vertex3}(1).
Since $\Gamma$ is not a cycle, the former is impossible.
Thus there is a good vertex of degree two.

If $\Gamma$ has only one pinched vertex,
then we see that there are five possibilities for $\Gamma$ as shown in Figure \ref{vertex3}(2)-(6).
Hence $\Gamma$ has a good vertex of degree at most $3$.

Finally, assume that $\partial\Gamma$ consists of two cycles.
Since there is no interior vertex, one cycle contains two vertices, and the other contains one vertex. 
Notice that any vertex is good.
By an Euler calculation, $\Gamma$ has at most $6$ edges.
Then $\Gamma$ contains (7) of Figure \ref{vertex3} as its subgraph.
If $z$ has degree $4$, then there are two edges connecting $z$ with $y$.
Then $x$ has degree $3$.
Thus $\Gamma$ has a good vertex of degree at most $3$.
\end{proof}

\section{The generic case}\label{generic}

In this section, we assume that $s\ge 4$ and $t\ge 4$.

\begin{lemma}
Any vertex of the reduced graph $\overline{G}_T$ has degree at least $5$.
\end{lemma}

\begin{proof}
Let $v$ be a vertex of $\overline{G}_T$.
If $\mathrm{deg}(v)\le 4$ in $\overline{G}_T$, then $\mathrm{deg}(v)\le 4s < 5s$ in $G_T$ by Lemma \ref{G}.
\end{proof}

\subsection{Some vertex of $\overline{G}_T$ has degree $5$}\label{case1}

In this subsection, we consider the case where some vertex of $\overline{G}_T$ has degree $5$, and show that
the case is impossible.
Let $v_i$ be such a vertex.

\begin{lemma}
In $G_T$, $v_i$ is incident to exactly five families of parallel negative edges, each of which contains $s$ edges.
\end{lemma}

\begin{proof}
This immediately follows from Lemma \ref{G}.
\end{proof}

Thus all $i$-edges in $G_S$ are positive by the parity rule, and there are five positive $i$-edges at each vertex of $G_S$.
Recall that $\overline{G}_S^+$ is the subgraph of $\overline{G}_S$ consisting all vertices and all positive edges.

\begin{lemma}\label{atleast5}
Any vertex of $\overline{G}_S^+$ has degree at least $5$.
\end{lemma}

\begin{proof}
This is because two positive $i$-edges at any vertex cannot be parallel by Lemma \ref{G}.
\end{proof}

From Proposition \ref{DandAsupp} and Lemma \ref{atleast5}, 
$\overline{G}_S^+$ has an interior vertex of degree at most $5$.
But the next lemma shows that this is impossible.

\begin{lemma}
$\overline{G}_S^+$ has no interior vertex of degree at most $5$.
\end{lemma}

\begin{proof}
Let $u$ be an interior vertex of $\overline{G}_S^+$ of degree at most $5$.
By Lemma \ref{atleast5}, $u$ has exactly degree $5$ in $\overline{G}_S^+$.
Since each family of parallel positive edges contains at most $t/2+1$ edges by Lemma \ref{G}(1), 
$u$ has degree at most $5(t/2+1)$ in $G_S$.
Hence $5(t/2+1)\ge \mathrm{deg}(u)=5t$, and then $t\le 2$, a contradiction.
\end{proof}

\subsection{Each vertex of $\overline{G}_T$ has degree $6$}

By the previous subsection \ref{case1}, 
we know that each vertex of $\overline{G}_T$ has degree at least $6$.
Then an easy Euler characteristic argument shows that each vertex of $\overline{G}_T$ has degree exactly $6$.
(See \cite[Claim 3.2]{BZ2}.)

\begin{lemma}\label{S-label}
If $\overline{G}_S^+$ has an interior vertex, then
$G_S$ has an $S$-cycle with label $j$ for any label $j$.
\end{lemma}

\begin{proof}
Let $u_i$ be an interior vertex of $\overline{G}_S^+$.
Then only positive edges are incident to $u_i$ in $G_S$.
By the parity rule, all $i$-edges in $G_T$ are negative.

There are five negative $i$-edges at the $j$th vertex $v_j$ of $G_T$, and any two of them are not parallel by Lemma \ref{G}(2).
Thus $v_j$ is incident to at most one positive edge in $\overline{G}_T$.
Hence $v_j$ is incident to at least $5s-(s/2+1)=9s/2-1$ negative edges.
In $G_S$, this means that there are at least $9s/2-1$ positive $j$-edges.

From an Euler characteristic calculation, $\overline{G}_S$ has at most $3s$ edges.
Since $9s/2-1>3s$, there are two positive $j$-edges which are parallel in $G_S$.
Then they form an $S$-cycle with $j$ as a label by Lemmas \ref{Scharlemann}(1) and \ref{G}(1).
\end{proof}

\begin{proposition}\label{nointvertex}
$\overline{G}_S^+$ cannot have an interior vertex.
\end{proposition}

\begin{proof}
Assume that $\overline{G}_S^+$ has an interior vertex.
By Lemma \ref{S-label}, any of the label set $\{1,2,\dots,t\}$ is a label of an $S$-cycle in $G_S$.
If $t=4$, then $G_S$ has two $S$-cycles with disjoint label pairs,
which is impossible by Lemma \ref{Scharlemann}(2).

Assume $t\ge 6$.
We may assume that $\{1,2\}$ is the label pair of an $S$-cycle of $G_S$.
Since $4$ is a label of an $S$-cycle, either $\{3,4\}$ or $\{4,5\}$ is the label pair of an $S$-cycle.
By Lemma \ref{Scharlemann}(4), it must be $\{3,4\}$.
Similarly, we can conclude that $\{5,6\}$ is the label pair of an $S$-cycle.
Thus there are three $S$-cycles with mutually disjoint label pairs, which is impossible by Lemma \ref{Scharlemann}(3).
\end{proof}

\begin{lemma}\label{key}
Let $u_i$ be a vertex of $\overline{G}_S^+$.
Suppose that some label $j$ appears $k$ times among negative edge endpoints of $u_i$ in $G_S$.
Then $k\le 4$.
Furthermore, if $k=4$ then $s=4$, and if $k=3$ then $s=4$ or $6$.
\end{lemma}

\begin{proof}
By the parity rule, there are $k$ positive $i$-edges at the vertex $v_j$ in $G_T$.
No two of them are parallel by Lemma \ref{G}(1).
Hence $k(s/2+1)+(6-k)s\ge \mathrm{deg}(v_j)=5s$.

If $k=5$, then $s\le 10/3$, a contradiction.
Thus we have $k\le 4$.
The others immediately follow from the inequality.
\end{proof}

\begin{lemma}\label{isolated}
$\overline{G}_S^+$ cannot have a vertex of degree at most one.
\end{lemma}

\begin{proof}
Assume that $u$ is a vertex of $\overline{G}_S^+$ of degree at most one.
Then there are at most $t/2+1$ positive edge endpoints at $u$ in $G_S$.
Hence at least $5t-(t/2+1)=9t/2-1$ negative edges are incident to $u$ successively.
Since $9t/2-1>4t$, some label appears five times among negative edge endpoints of $u$.
This is impossible by Lemma \ref{key}.
\end{proof}

\begin{lemma}\label{degree2}
If $\overline{G}_S^+$ has a vertex $u_i$ of degree two, then $s=4$, and $G_T$ has an $S$-cycle with $i$ as a label.
\end{lemma}

\begin{proof}
Since there are at most $2(t/2+1)=t+2$ positive edge endpoints at $u_i$ in $G_S$,
$u_i$ has at least $4t-2$ negative edge endpoints.
Then $4t-2>3t$, and hence some label $j$ appears at least $4$ times among negative edge endpoints of $u_i$.
By Lemma \ref{key}, $s=4$.

Also, $G_T$ has $4t-2$ positive $i$-edges.
Since $\overline{G}_T$ has at most $3t$ edges (as seen by an Euler characteristic calculation),
some positive $i$-edges are parallel in $G_T$.
Thus $G_T$ has an $S$-cycle with $i$ as a label by Lemmas \ref{Scharlemann}(1) and \ref{G}(1).
\end{proof}

\begin{proposition}\label{noD-supp}
$\overline{G}_S^+$ has no component with a disk support.
\end{proposition}

\begin{proof}
Assume not.
Choose an extremal component $\Gamma$ with a disk support.
Then $\Gamma$ has a good vertex of degree at most two by Lemma \ref{Dsupp}.
Hence $s=4$ by Lemma \ref{degree2}.
Thus $\Gamma$ has at most two vertices.
But this is impossible by Proposition \ref{nointvertex} and Lemma \ref{isolated}.
\end{proof}

By Lemma \ref{DorA} and Proposition \ref{noD-supp}, 
any component of $\overline{G}_S^+$ has an annulus support and
there are at least two components.

\begin{lemma}\label{no2}
$\overline{G}_S^+$ has no cycle component.
\end{lemma}

\begin{proof}
Assume that $\overline{G}_S^+$ has a cycle component $\Gamma$.
By Lemma \ref{degree2}, $s=4$.
Hence $\Gamma$ contains at most two vertices.
Recall that $G_S$ has vertices $u_1,u_2,u_3,u_4$, where
$u_i$ and $u_j$ are parallel if and only if $i\equiv j\pmod{2}$.

First, assume that $\Gamma$ contains only one vertex.
Then we can assume that $\overline{G}_S^+$ has two loop components (with annulus supports) based on $u_1$ and $u_3$.
By Lemma \ref{degree2}, $G_T$ has $S$-cycles with labels $1$ and $3$, respectively.
Moreover, we can assume that their label pairs are $\{1,2\}$ and $\{2,3\}$ by Lemma \ref{Scharlemann}(2).

Then $\overline{G}_S$ contains a subgraph as shown in Figure \ref{nocycle} by Lemma \ref{Schsupp}.
By \cite[Lemma 1.9]{Wu}, $u_4$ must lie in the disk region $D$ as indicated in Figure \ref{nocycle}.
Then  $\overline{G}_S^+$ has a component containing both of $u_2$ and $u_4$.
Hence $u_4$ has degree two in $\overline{G}_S^+$.
By Lemma \ref{degree2} again, $G_T$ has an $S$-cycle with $4$ as a label, that is, an $S$-cycle with label pair either $\{3,4\}$ or $\{4,1\}$.
In either case, this contradicts Lemma \ref{Scharlemann}(2).

\begin{figure}[tb]
\includegraphics*[scale=0.3]{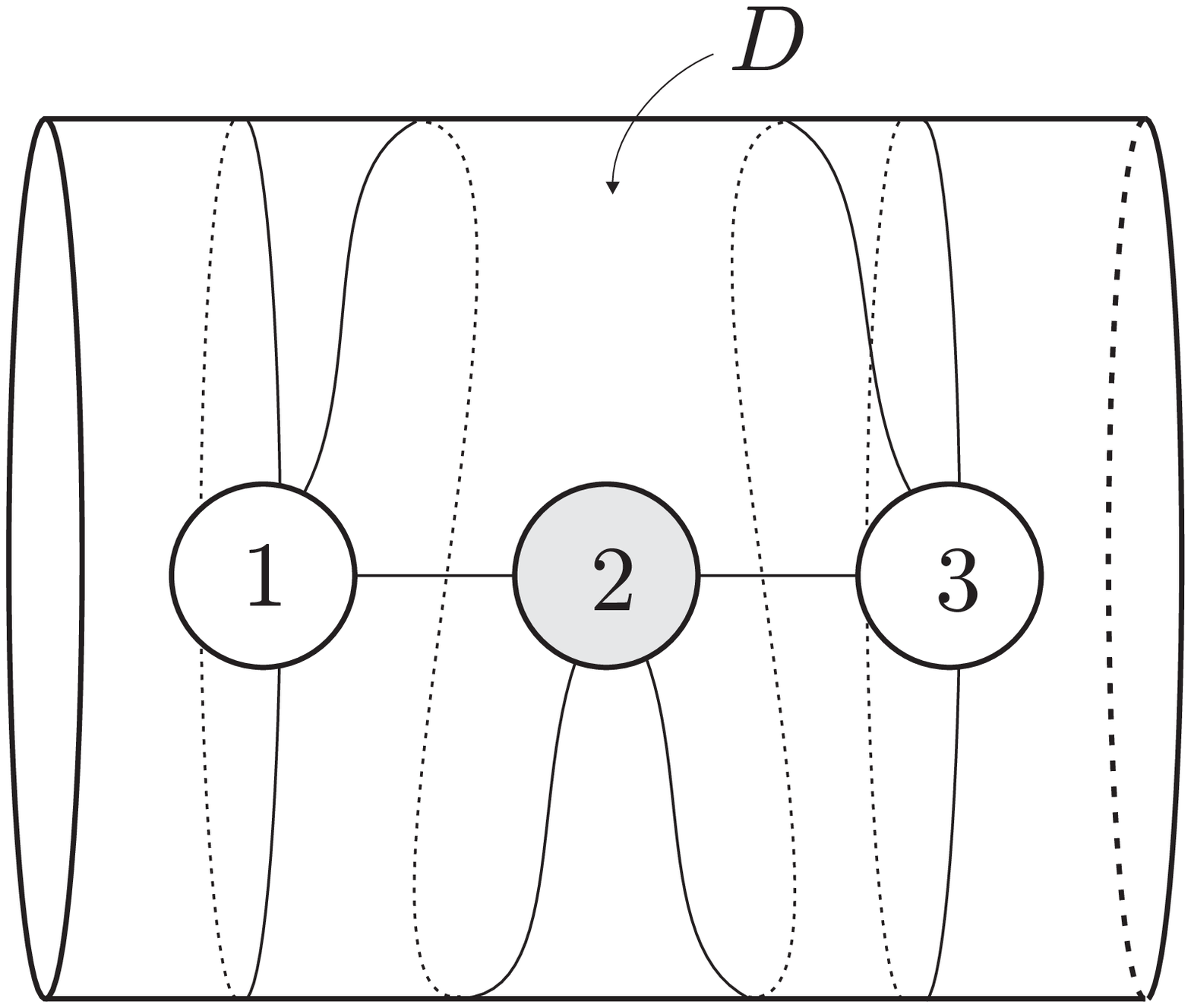}
\caption{}\label{nocycle}
\end{figure}

Next, assume that $\Gamma$ contains just two vertices, $u_1$ and $u_3$, say.
Again, $G_T$ has two $S$-cycles $\rho_1$ and $\rho_2$ with label pairs $\{1,2\}$ and $\{2,3\}$, respectively.
Also, $\overline{G}_S^+$ has another component $\Lambda$
containing $u_2$ and $u_4$, otherwise $u_4$ has degree two, which leads to a contradiction as above.
In fact, $\Lambda$ has the form of either (2), (3) or (4) in Figure \ref{vertex2}, where $u_4$ has degree at least $3$.
By Lemma \ref{Schsupp}, the edges of $\rho_i$ form an essential loop on $\widehat{S}$.
Then we cannot place the edges of two $S$-cycles $\rho_1$ and $\rho_2$ to satisfy this condition simultaneously.
\end{proof}

Let $\Gamma$ be an extremal component of $\overline{G}_S^+$.
It has an annulus support, and it is not a cycle by Lemma \ref{no2}.
Therefore $\Gamma$ has a good vertex of degree at most $4$ by Lemma \ref{Asupp}.
Let $u$ be such a vertex.

\begin{lemma}
$s=4$ or $6$.
\end{lemma}

\begin{proof}
There are at most $4(t/2+1)=2t+4$ positive edge endpoints at $u$ in $G_S$ by Lemma \ref{G}(1).
Thus $u$ has at least $3t-4$ negative edge endpoints.  If $t>4$, then $3t-4>2t$.
If $t=4$, then there are at most $3(t/2+1)+t/2=2t+3$ positive edge endpoints at $u$ in $G_S$ by Lemma \ref{consec}.
Thus $u$ has at least $3t-3$ negative edge endpoints, and note $3t-3>2t$.
Hence, in either case, some label appears at least three times among negative edge endpoints of $u$.
Then $s=4$ or $6$ by Lemma \ref{key}.
\end{proof}

Now, we divide the cases.


\subsubsection*{Case 1. $s=6$}

\begin{lemma}
$\overline{G}_S^+$ consists of two components, each of which has three vertices.
\end{lemma}

\begin{proof}
Recall that $\overline{G}_S^+$ consists of at least two components, each of which has an annulus support.
Since there is no cycle in $\overline{G}_S^+$ by Lemma \ref{no2}, each component must contain three vertices.
\end{proof}

\begin{proposition}\label{s6}
$s=6$ is impossible.
\end{proposition}

\begin{proof}
By Lemmas \ref{3vertex}, \ref{isolated}, \ref{degree2} and \ref{no2}, $\overline{G}_S^+$ has a good vertex $u_i$ of degree $3$.
Assume $t\ge 8$.
By Lemma \ref{G}(1), there are at most $3(t/2+1)$ positive edge endpoints at $u_i$ in $G_S$.
Thus $u_i$ has at least $7t/2-3$ negative edge endpoints.
Since $7t/2-3>3t$, some label appears four times among negative edge endpoints of $u_i$.
Then $s=4$ by Lemma \ref{key}, which is a contradiction.

Assume $t=6$.
It suffices to consider the case where $u_i$ is incident to three families of $4$ parallel positive edges.
(Otherwise, there are more than $18(=3t)$ negative edge endpoints at $u_i$, and then some label appears four times there.)
Then $18$ negative edges are incident to $u_i$ successively in $G_S$.
Thus any label $j$ appears exactly three times there.
In fact, three occurrences of the label $j$ are consecutive among five occurrences of $j$ at $u_i$.
In $G_T$, there are three positive $i$-edges at $v_j$, whose endpoints with label $i$ are
consecutive at $v_j$ among the five occurrences of label $i$ by Lemma \ref{jumping}.
Since no two of the $i$-edges are parallel, $v_j$ is incident to three families of parallel positive edges containing $i$-edges, which
are consecutive.
Hence $v_j$ has at least $13$ positive edge endpoints.
Thus $v_j$ is incident to at least $4$ families of parallel positive edges by Lemma \ref{G}(1). 
Then $v_j$ is incident to at most two families of parallel negative edges.
But this implies that $v_j$ has at most $4\cdot 4+6\cdot 2=28$ edge endpoints, which contradicts
that it has degree $30$.

Next assume $t=4$.
In $G_S$, $u_i$ has at most $9$ positive edge endpoints.
Hence there are at least $11$ negative edges there.
Thus some label appears three times among negative edge endpoints of $u_i$.
A similar argument to the case $t=6$ above leads to a contradiction.
We have thus shown that the case $s=6$ is impossible.
\end{proof}


\subsubsection*{Case 2. $s=4$}

By Lemmas \ref{isolated} and \ref{no2},
$\overline{G}_S^+$ consists of two connected components, each of which has the form of Figure \ref{vertex2}(2), (3) or (4)
by Lemma \ref{2vertex}.

\begin{lemma}\label{vertex3-1}
$\overline{G}_S^+$ does not have a component of the form as in Figure \ref{vertex2}\,\textup{(2)}.
\end{lemma}

\begin{proof}
Let $\Gamma$ be a component of $\overline{G}_S^+$ as in Figure \ref{vertex2}(2), and
let $u_i$ be the good vertex of degree two in $\Gamma$.
As in the proof of Lemma \ref{degree2}, some label $j$ appears four times among negative edge endpoints of $u_i$.
In $G_T$, there are four positive $i$-edges at $v_j$.
No two of them are parallel.
Thus $v_j$ is incident to four families of parallel positive edges, each of which contains an $i$-edge.
Then we see that those four families contain $3$ edges respectively, and that $v_j$ is incident to two
families of $4$ parallel negative edges.
By Lemma \ref{jumping}, the families of positive edge are consecutive.
But this contradicts Lemma \ref{consec}.
\end{proof}

\begin{lemma}\label{vertex3-2}
$\overline{G}_S^+$ does not have a component of the form as in Figure \ref{vertex2}\,\textup{(3)}.
\end{lemma}

\begin{figure}[tb]
\includegraphics*[scale=0.5]{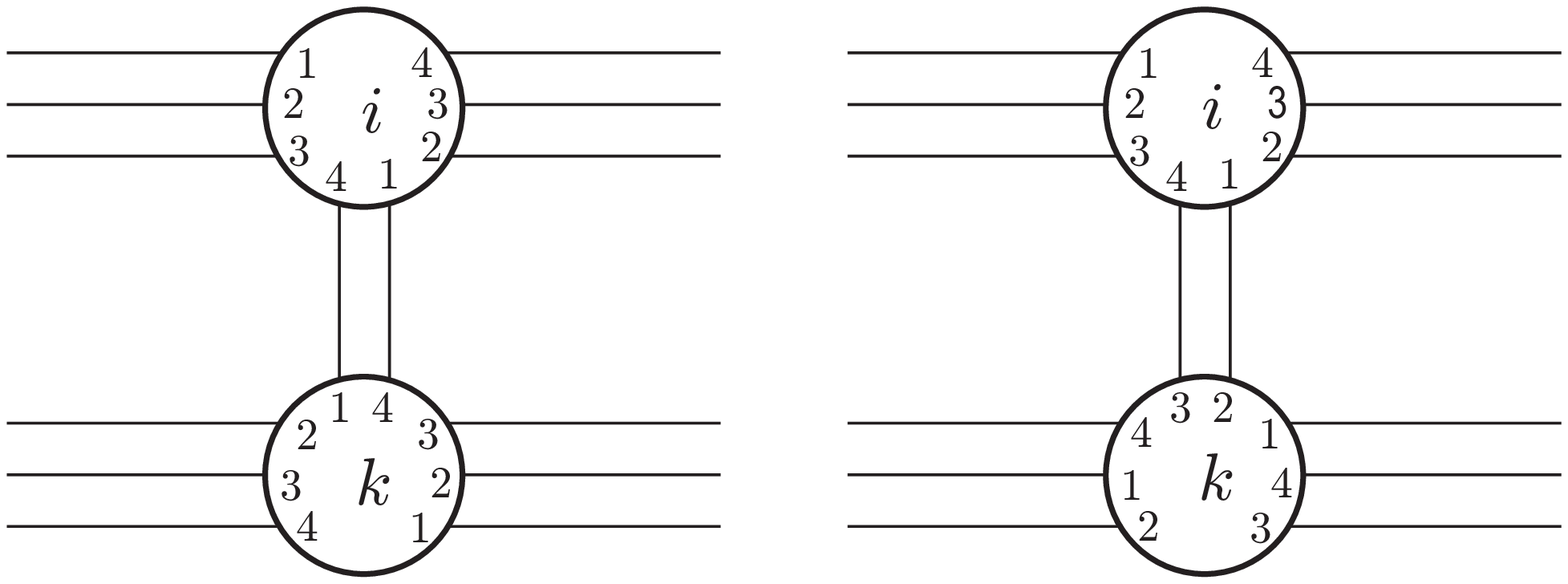}
\caption{}\label{s4case2}
\end{figure}

\begin{proof}
Let $\Gamma$ be such a component.
Then $\Gamma$ has a good vertex $u_i$ of degree $3$.
Assume $t\ge 8$.
As in the first paragraph of the proof of Proposition \ref{s6}, 
some label $j$ appear four times among negative edge endpoints of $u_i$.
Hence $v_j$ has $4$ positive $i$-edges, which are not mutually parallel.
Since $v_j$ has degree $20$ in $G_T$, $v_j$ is incident to four families of $3$ parallel positive edges and
two families of $4$ parallel negative edges.
Then the four families of parallel positive edges are consecutive by Lemma \ref{jumping}.
But this contradicts Lemma \ref{consec}.

Assume $t=6$.
It suffices to consider the case where $u_i$ is incident to three families of $4$ parallel positive edges.
(Otherwise, $u_i$ has more than $18(=3t)$ negative edge endpoints, and then
some label appears four times there.)
Then $u_i$ is incident to $4$ loops and a family of $4$ parallel positive non-loop edges.
We can assume that the labels at $u_i$ are as shown in Figure \ref{s4case1}.
Let $u_k$ be another vertex of $\Gamma$.
Then the situation at $u_k$ is the same as $u_i$.
Hence $\Gamma$ has three $S$-cycles with disjoint label pairs, which is impossible by Lemma \ref{Scharlemann}(3).

\begin{figure}[tb]
\includegraphics*[scale=0.5]{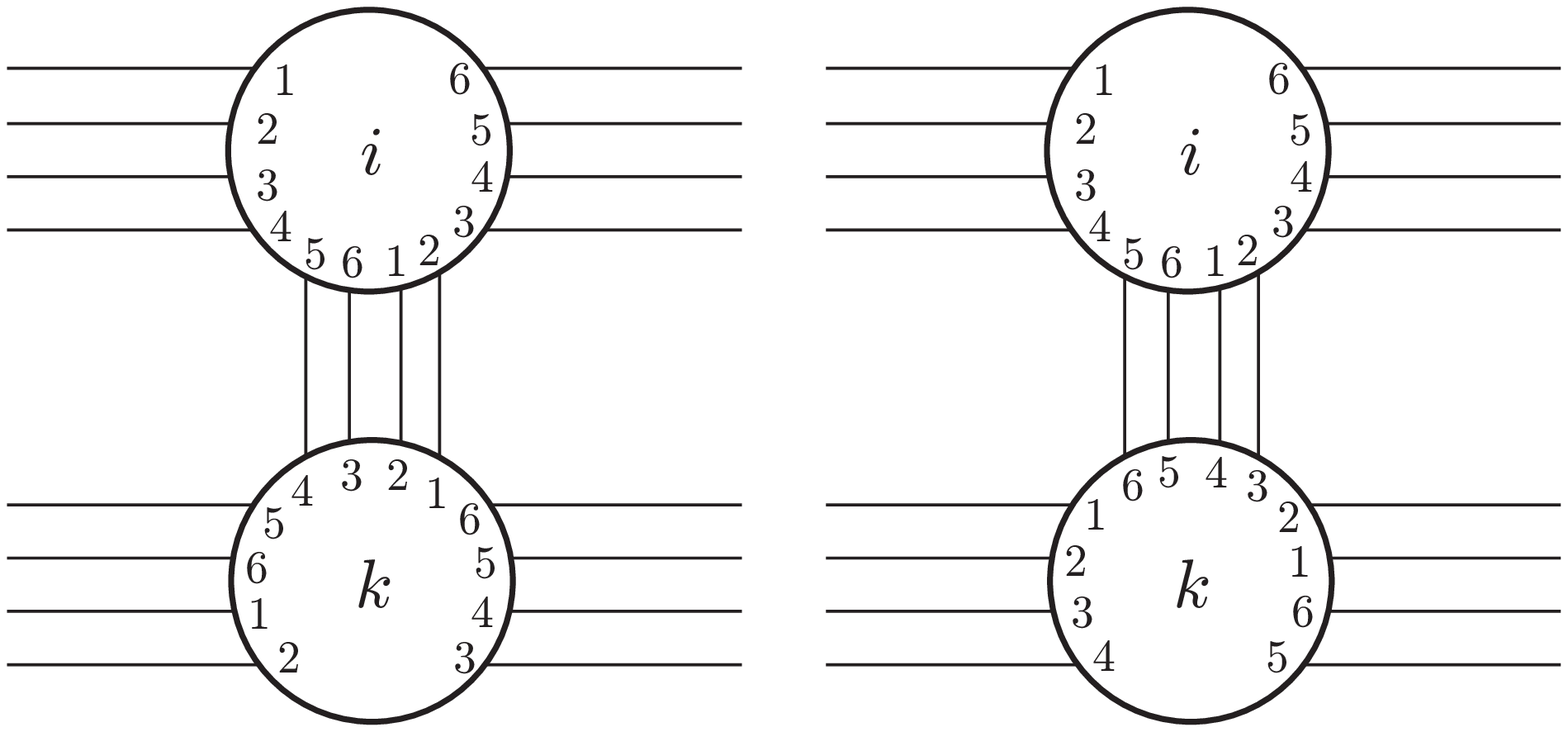}
\caption{}\label{s4case1}
\end{figure}

Assume $t=4$.
If $u_i$ is incident to more than $12$ negative edges, then
some label appears four times among negative edge endpoints of $u_i$.
This leads to a contradiction as above.
Thus $u_i$ has at most $12$ negative edges, and then there are $8$ or $9$ positive edge endpoints.
If there are $9$ positive edge endpoints at $u_i$, three loops and a family of three positive non-loop edges are incident to $u_i$.
But this contradicts the parity rule.
Hence $u_i$ has exactly $12$ negative edge endpoints and $8$ positive edge endpoints.
The parity rule implies that there are three loops and two non-loop edges.
We can assume that the labels at $u_i$ as shown in Figure \ref{s4case2}.
Then $\Gamma$ has two $S$-cycles with disjoint label pairs as in Figure \ref{s4case2}, which is impossible by Lemma \ref{Scharlemann}(2).
\end{proof}

\begin{proposition}
$s=4$ is impossible.
\end{proposition}

\begin{proof}
By Lemmas \ref{vertex3-1} and \ref{vertex3-2}, $\overline{G}_S^+$ consists of two components of the form as in Figure \ref{vertex2}(4).
Let $u$ be any vertex of $G_S$.
Since at most two families of parallel negative edges are incident to $u$, there are at most $2t$ negative edge endpoints at $u$.
Thus $u$ has at least $3t$ positive edge endpoints.
Then $4(t/2+1)\ge 3t$, and so $t=4$.
Hence $u$ has three loops and two families of three parallel positive edges.
Then $G_S$ has two $S$-cycles with disjoint label pairs, which is impossible by Lemma \ref{Scharlemann}(2).
Thus we have shown that $s=4$ is impossible.
\end{proof}

\section{The case that $s=2$ and $t\ge 4$}\label{s2t4}

In this section, we assume $s=2$ and $t\ge 4$, but all arguments can apply to
the case that $t=2$ and $s\ge 4$.

\begin{figure}[tb]
\includegraphics*[scale=0.35]{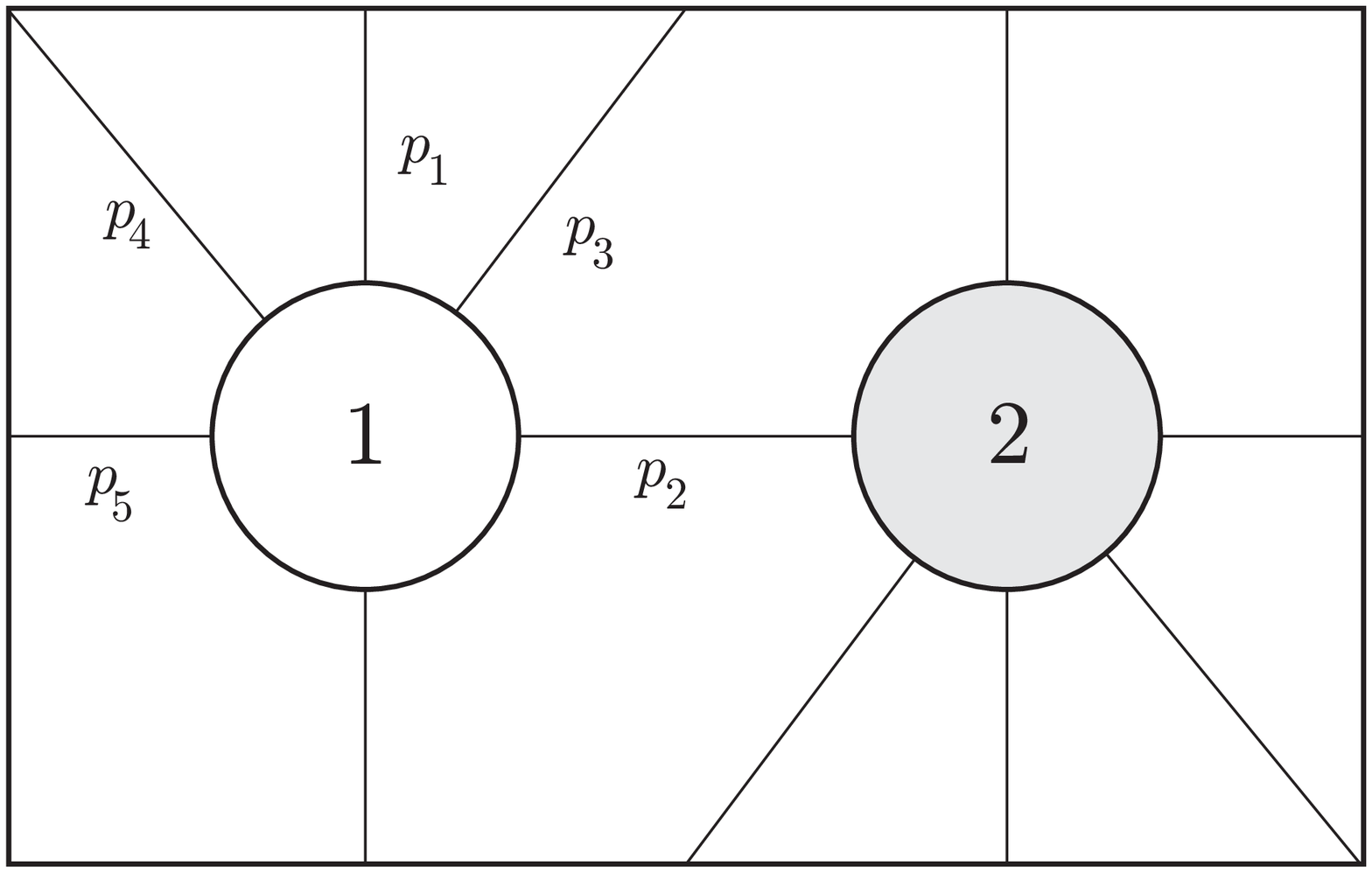}
\caption{}\label{t2}
\end{figure}

The reduced graph $\overline{G}_S$ is a subgraph of the graph as shown in Figure \ref{t2} \cite[Lemma 5.2]{G},
where the sides of the rectangle are identified to form $\widehat{S}$ in the usual way.
Here, $p_i$ indicates the number of edges in the family of parallel edges.
Recall that $p_1\le t/2+1$ and $p_i\le t$ for $i=2,3,4,5$ by Lemma \ref{G}.

\begin{lemma}
In $\overline{G}_S$, $u_1$ and $u_2$ have degree $6$.
Moreover, $p_1=t/2$ or $t/2+1$.
\end{lemma}

\begin{proof}
Since $u_i$ has degree $5t$ in $G_S$, the first follows immediately.
Also, $u_i$ has at least $t$ positive edge endpoints.  Thus the second follows.
\end{proof}

We distinguish two cases.

\subsection*{Case 1. $p_1=t/2$}

\begin{figure}[tb]
\includegraphics*[scale=0.4]{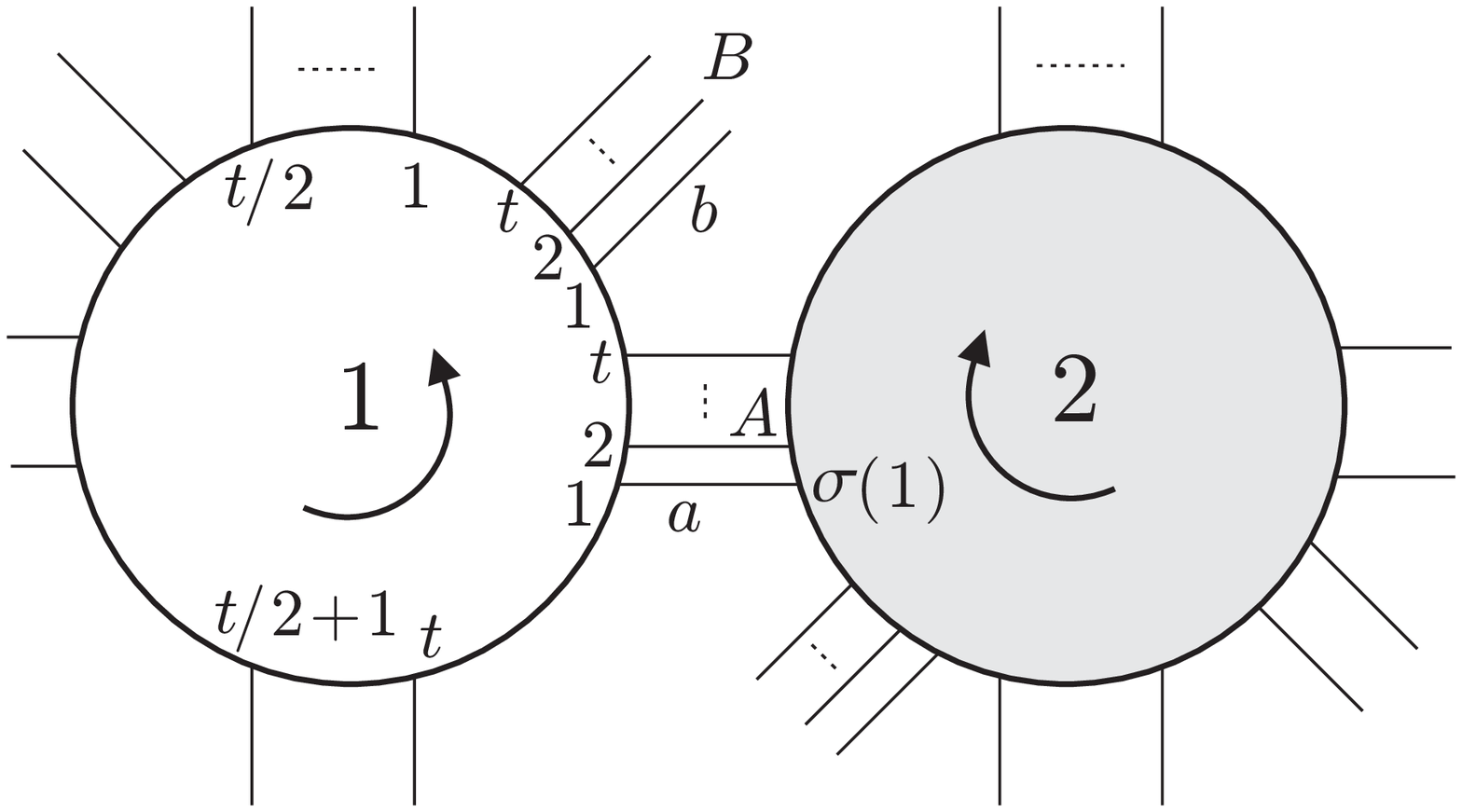}
\caption{}\label{t2case1}
\end{figure}

In this case, $p_i=t$ for $i=2,3,4,5$.
Let $A$ and $B$ be the families of $p_2$ and $p_3$ parallel negative edges in $G_S$, respectively.
We can assume that the labels are as in Figure \ref{t2case1}.
Let $\sigma$ be the associated permutation to $A$ such that
an edge in $A$ has label $j$ at $u_1$ and label $\sigma(j)$ at $u_2$.
The edges of $A$ form disjoint cycles in $G_T$ according to the orbits of $\sigma$, and such a cycle is essential on $\widehat{T}$ \cite[Lemma 2.3]{G}.
By the parity rule, each cycle contains only the vertices of the same sign.
Hence there are at least two such cycles.
Let $L$ be the cycle corresponding to the orbit of $\sigma$ containing the label $1$.

Note that the four families of negative edges in $G_S$ define the same permutation $\sigma$.

\begin{lemma}\label{notid}
$\sigma$ is not the identity.
\end{lemma}

\begin{figure}[tb]
\includegraphics*[scale=0.4]{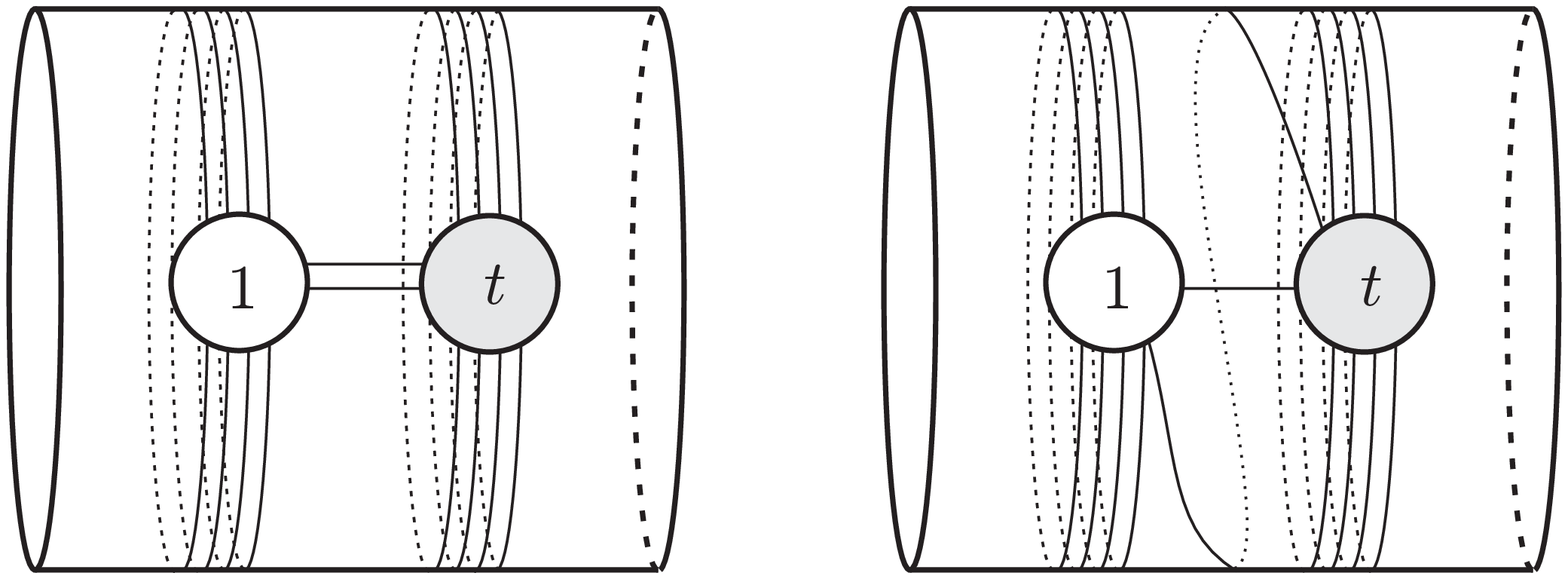}
\caption{}\label{s2case1id}
\end{figure}

\begin{proof}
Assume that $\sigma$ is the identity.
Then each family of parallel negative edges in $G_S$ contains a $\{j,j\}$-edge for $j=1,2,\dots,t$.
Let $G(1,t)$ be the subgraph of $G_T$ spanned by the vertices $v_1$ and $v_t$.
Then $G(1,t)$ has an annulus support on $\widehat{T}$, since $t\ge 4$.
Hence there are two possibilities for $G(1,t)$ as shown in Figure \ref{s2case1id}.

But a jumping number argument will eliminate both configurations as follows.
Let $a$ be the $\{1,1\}$-edge in $A$, and let $a_i$ be its endpoint at $u_i$ for $i=1,2$.
There are two positive $\{1,t\}$-loops $e$ based on $u_1$ and $f$ based on $u_2$ in $G_S$.
Let $e_1$ and $f_1$ be their endpoints with label $1$.
Around $u_1$, $a_1$ and $e_1$ are not successive among five occurrences of label $1$, but
$a_2$ and $f_1$ are successive among five occurrences of label $1$ around $u_2$.
By Lemma \ref{jumping}, $a_1$ and $e_1$ are not successive among five occurrences of label $1$ around $v_1$, but
$a_2$ and $f_1$ are successive among five occurrences of label $2$ around $v_1$.
But this is not satisfied in both configurations of $G(1,t)$.
\end{proof}

\begin{lemma}\label{involution}
$\sigma^2$ is the identity.  In particular, each orbit of $\sigma$ has length two, and $\sigma(1)=t/2+1$.
\end{lemma}

\begin{proof}
Let $a$ ($b$ resp.) be the edge of $A$ ($B$ resp.) with label $1$ at $u_1$.
Then $L\cup b$ is contained in an annulus on $\widehat{T}$.
There are two possibilities for $L\cup b$ as shown in Figure \ref{s2case1inv}, where we put $r=\sigma(1)$.
Note that $a$ and $b$ have label $1$ at $v_1$.

\begin{figure}[tb]
\includegraphics*[scale=0.35]{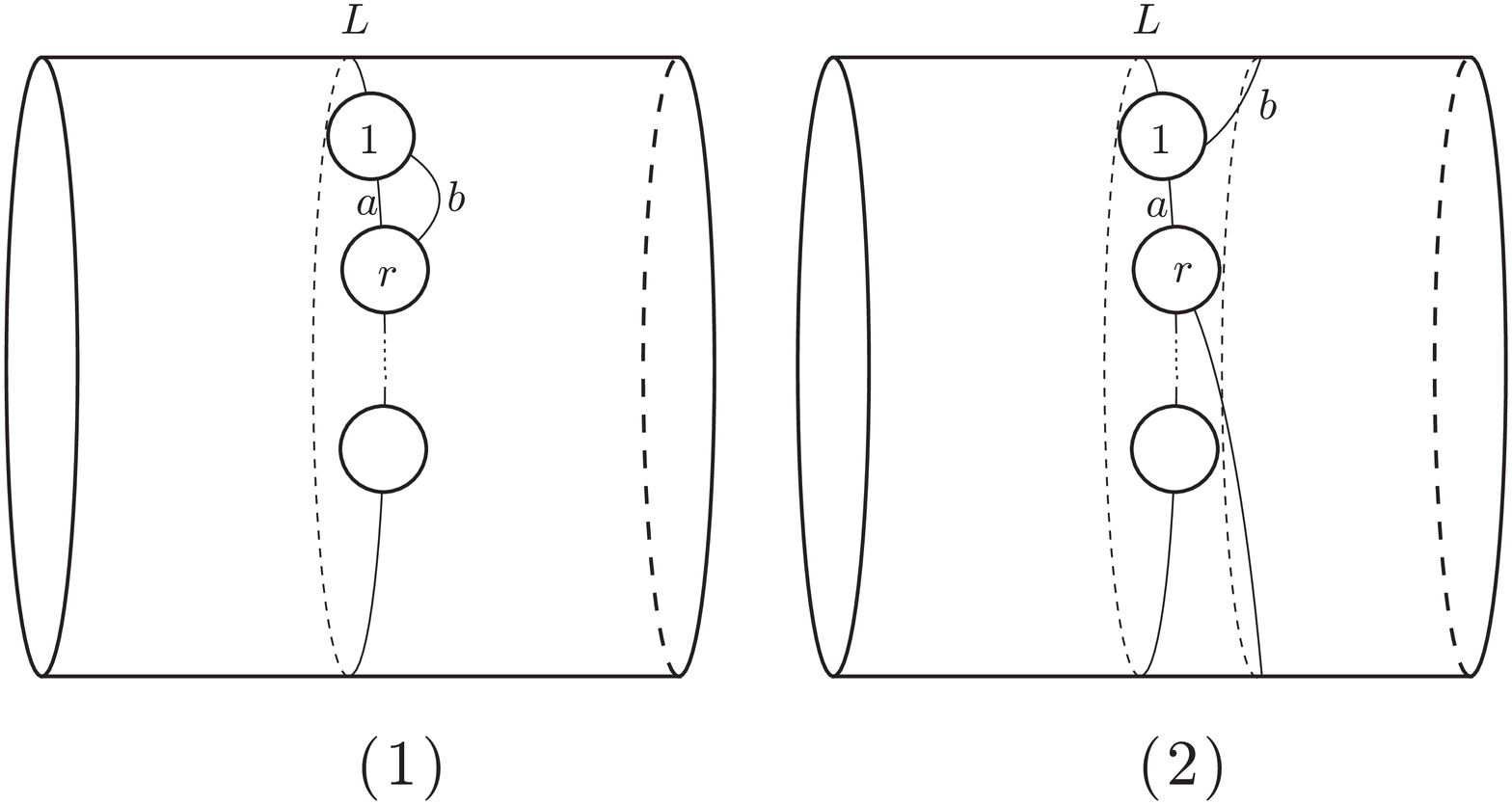}
\caption{}\label{s2case1inv}
\end{figure}

For Figure \ref{s2case1inv}(1), there is another edge $e$ between $a$ and $b$.
Then $e$ is a negative $\{1,r\}$-edge in $G_S$ with label $r$ at $u_1$ and label $1$ at $u_2$.
Although $e$ need not be in $A$, this implies $\sigma(r)=1$, because any family of negative edges corresponds to the same permutation $\sigma$.
Hence $\sigma^2$ is the identity.

For Figure \ref{s2case1inv}(2), suppose that $\sigma^2$ is not the identity.
Then $L$ contains at least three vertices.
Let $c$ ($d$, resp.) be the edge in $A$ ($B$, resp.) with label $r$ at $u_1$.
Of course, $c$ is contained in $L$.
Then $d$ and $b$ are on the same side of $L$, because the endpoints of $c$ and $d$ with label $r$ are successive around $u_1$.
Hence $d$ is parallel to $c$ in $G_T$.
Then there is another edge between them, which implies $\sigma^2$ is the identity as above.
This is a contradiction.
\end{proof}

\begin{lemma}\label{t=4}
$t=4$.
Furthermore, $G_T$ has a torus support.
\end{lemma}

\begin{proof}
By Lemma \ref{involution}, $G_T^+$ consists of $t/2$ components, and hence each component has an annulus support.
Let $G^+(1,t/2+1)$ be the component of $G_T^+$ containing the vertices $v_1$ and $v_{t/2+1}$.
Then it consists of $8$ edges, which are split into two families of $4$ parallel edges.
For, if a family contains $5$ edges, then some pair of edges is parallel in $G_S$, too.
This contradicts Lemma \ref{bothparallel}.
Similarly, let $G^+(t/2,t)$ be the component of $G_T^+$ containing $v_{t/2}$ and $v_t$.
Since $G_S$ has two $\{1,t\}$-loops and two $\{t/2,t/2+1\}$-loops (see Figure \ref{t2case1}),
the component $H$ of $G_T$ containing $G^+(1,t/2+1)$ and $G^+(t/2,t)$ has the form as shown in Figure \ref{s2case1tg4}, 
under the assumption $t>4$.

\begin{figure}[tb]
\includegraphics*[scale=0.35]{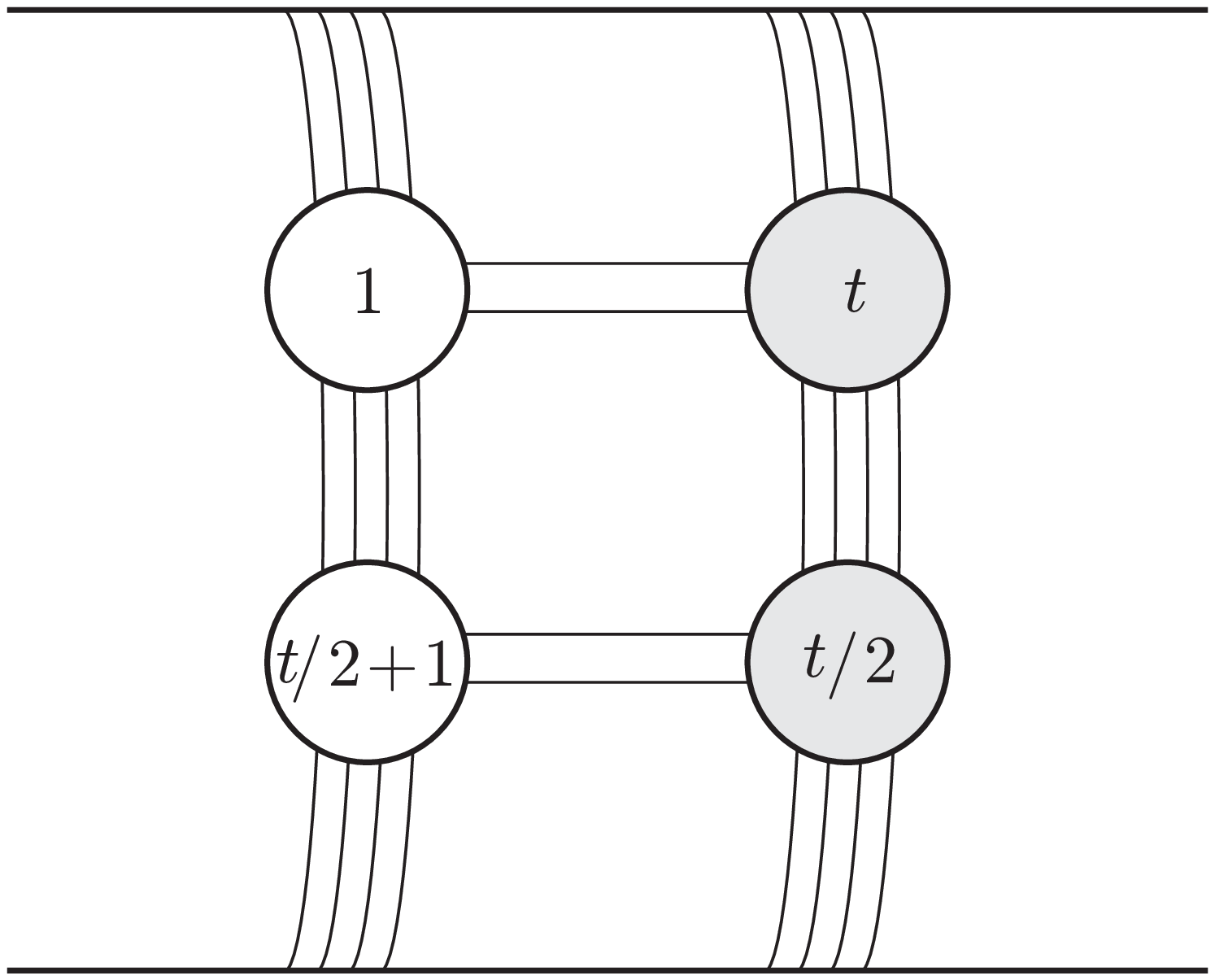}
\caption{}\label{s2case1tg4}
\end{figure}

But a jumping number argument will eliminate this configuration as before.
Look at the edge $a$ in the family of $A$ in $G_S$.
The endpoint of $a$ at $u_1$ is not adjacent to the endpoint of the $\{1,t\}$-loop with label $1$ among five occurrences of label $1$.
But the endpoint of $a$ at $u_2$ is adjacent to the endpoint of the $\{t/2,t/2+1\}$-loop with label $t/2+1$ among five
occurrences of label $t/2+1$.
Then we cannot locate the edge $a$ in $H$ to satisfy Lemma \ref{jumping}.
Hence $t=4$.

Then, $G_T^+=G^+(1,3)\cup G^+(2,4)$.
In this case, $H=G_T$.  If $H$ has an annulus support, then we have a contradiction as above.
Thus $G_T$ has a torus support.
\end{proof}

\begin{figure}[tb]
\includegraphics*[scale=0.4]{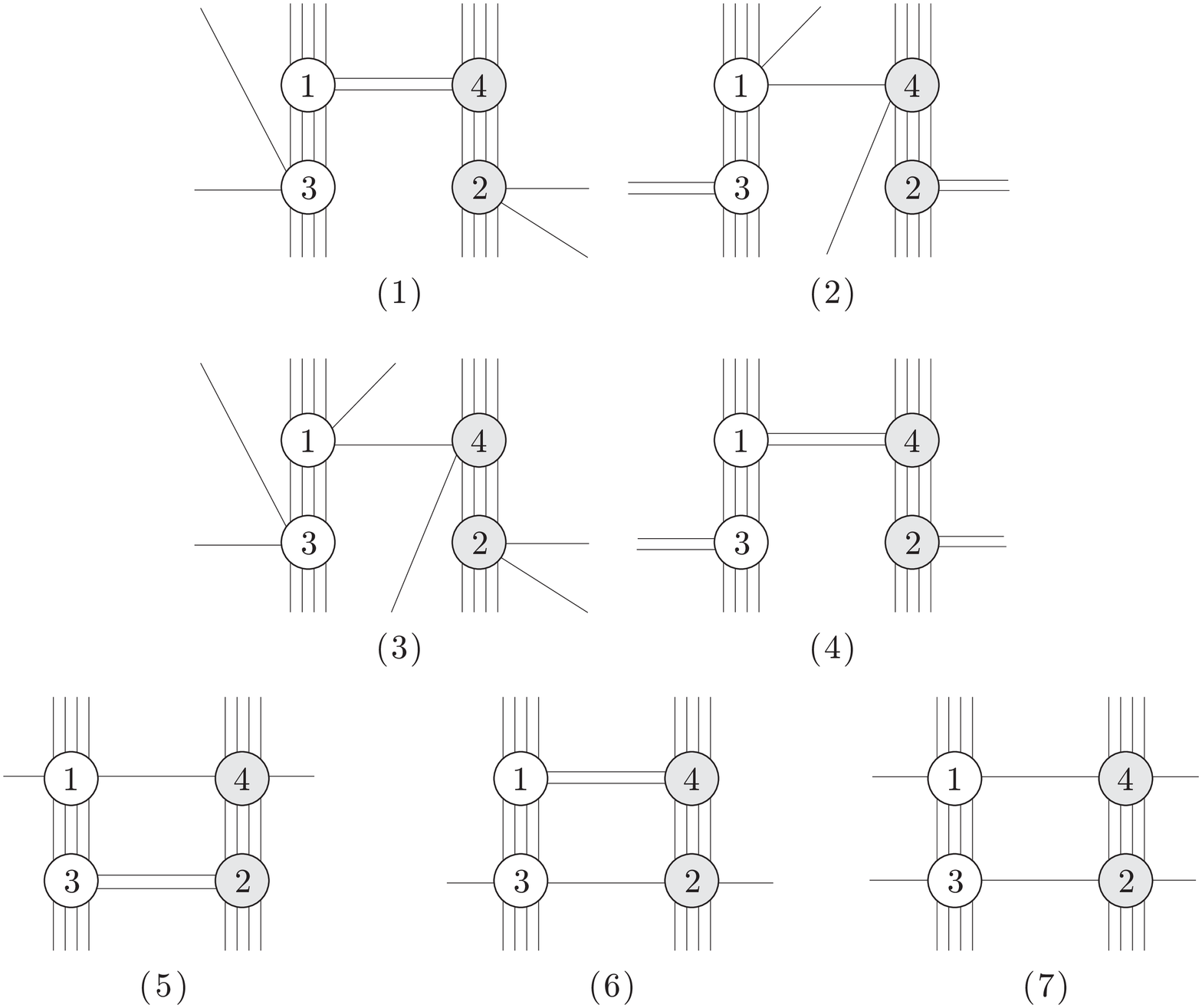}
\caption{}\label{s2case1t4}
\end{figure}

Thus $G_S$ is uniquely determined, and then there are seven possibilities for $G_T$ as shown in Figure \ref{s2case1t4}.
Clearly, (1), (2), (5) and (6) contradict the parity rule.

\begin{lemma}
\textup{(3)} of Figure \ref{s2case1t4} is impossible.
\end{lemma}

\begin{figure}[tb]
\includegraphics*[scale=0.63]{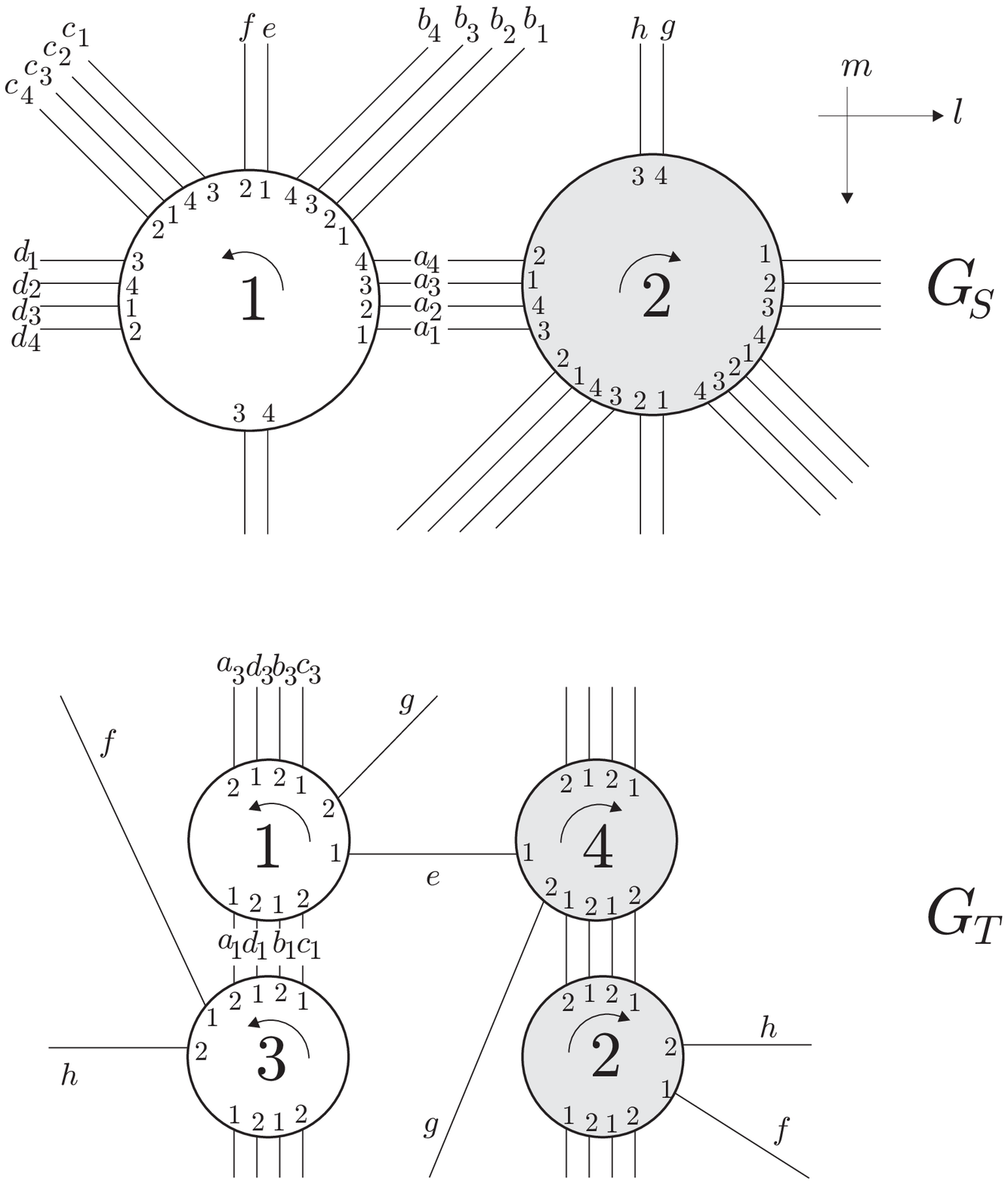}
\caption{}\label{worst1}
\end{figure}

\begin{proof}
The (partial) correspondence between the edges of $G_S$ and $G_T$ are determined as in Figure \ref{worst1} by using Lemma \ref{jumping}. 
Let $V_{12}$ be the part of $V_\beta$ between vertices $v_1$ and $v_2$ (disjoint from $v_3$ and $v_4$).
Ten edge endpoints at $v_1$ are connected to those at $v_2$ by disjoint arcs on the annulus $\mathrm{cl}(\partial V_{12}-v_1\cup v_2)$.
In particular, the consecutive endpoints of $e$, $g$ and $c_3$ at $v_1$ are connected to
the consecutive ones of $f$, $h$ and $c_4$ at $v_2$, respectively.
Also, the anticlockwise ordering of the former at $v_1$ must determine the clockwise ordering of the latter at $v_2$.
This contradicts Figure \ref{worst1}.
\end{proof}

\begin{lemma}
\textup{(7)} of Figure \ref{s2case1t4} is impossible.
\end{lemma}

\begin{proof}
We use the notation of the edges of $G_S$ in Figure \ref{worst1}.
At $u_1$, the endpoints of $e$ and $c_3$ are adjacent among five occurrences of label $1$. 
By Lemma \ref{jumping}, the endpoints of $e$ and $c_3$ are adjacent among five occurrences of label $1$ at $v_1$.
Then the endpoints of $h$ and $c_3$ are adjacent among five occurrences of label $2$ at $v_3$.
But this leads to a contradiction, because the endpoints of $h$ and $c_3$ are not adjacent among five
occurrences of label $3$ at $u_2$.
\end{proof}


\begin{lemma}\label{lem:worst2}
\textup{(4)} of Figure \ref{s2case1t4} is impossible.
\end{lemma}

\begin{proof}
The correspondence between the edges of $G_S$ and $G_T$ are determined as in Figure \ref{worst2} by using Lemma \ref{jumping},
where $G_S$ is the same as in Figure \ref{worst1}. 
To calculate $H_1(K(\alpha))$, we build up $K(\alpha)$ based on $\widehat{S}\cup V_\alpha$.

\begin{figure}[tb]
\includegraphics*[scale=0.63]{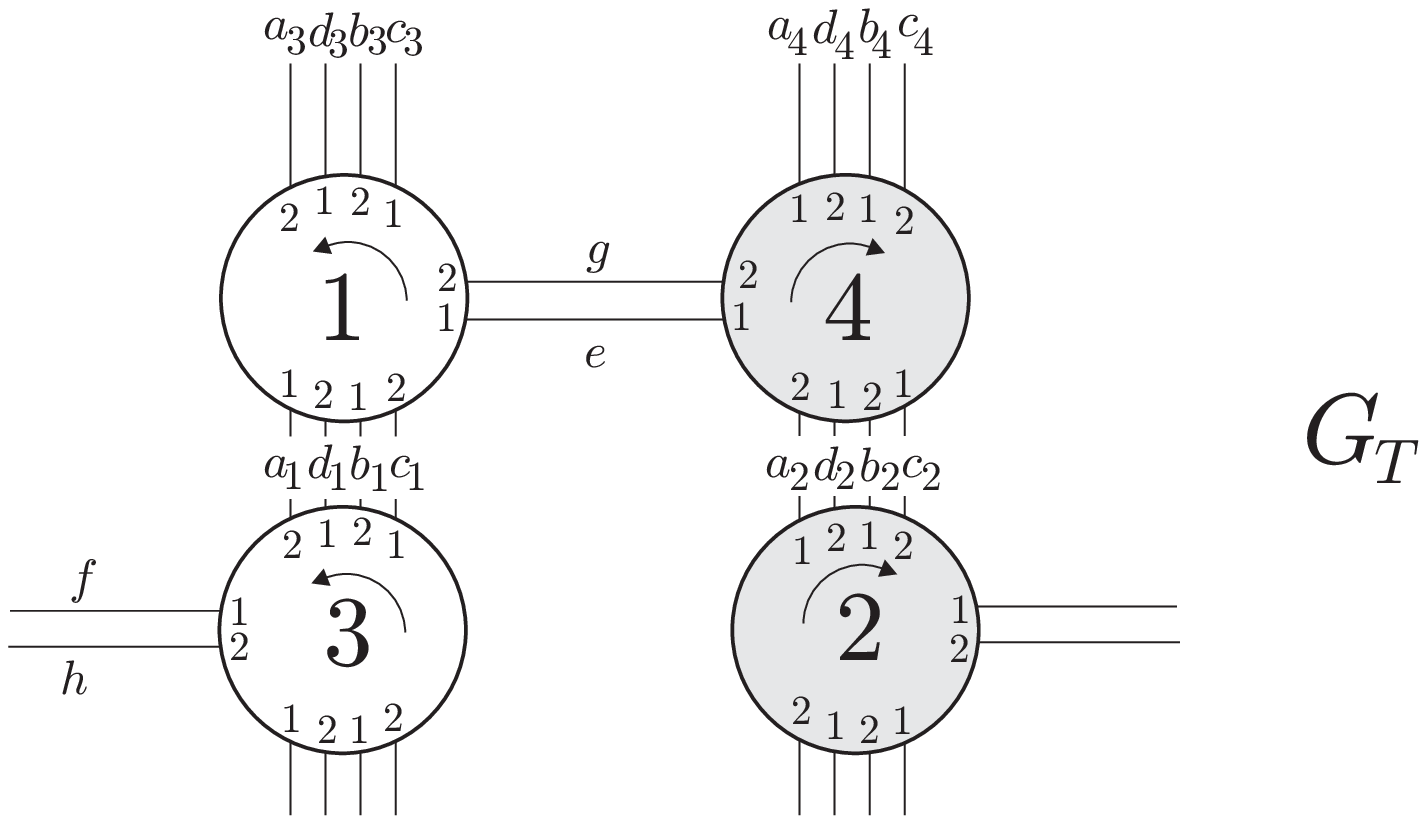}
\caption{}\label{worst2}
\end{figure}

Consider the bigons $D_1$ and $D_2$ in $G_T$ between $a_1$ and $d_1$, $e$ and $g$, respectively.
They are on the same side of $\widehat{S}$.
Let us call this side $\mathcal{B}$, the other side $\mathcal{W}$.
Thus $K(\alpha)=\mathcal{B}\cup \mathcal{W}$, and $\mathcal{B}\cap \mathcal{W}=\widehat{S}$.
Let $V_{12}=V_\alpha\cap \mathcal{B}$ and $V_{21}=V_\alpha\cap \mathcal{W}$.
Let $F$ be the genus two closed surface obtained from $\widehat{S}$ by tubing along $V_{12}$.
That is, $F=(\widehat{S}-u_1\cup u_2)\cup H$, where $H$ is the annulus $\mathrm{cl}(\partial V_{12}-u_1\cup u_2)$.
On $F$, $\partial D_1$ is non-separating, because it runs twice on $V_{12}$ in the same direction.
Hence surgering $F$ along $D_1$ gives a torus $\widehat{S}'$.
Furthermore, we see that $\partial D_2$ is non-separating on $\widehat{S}'$.
Hence $\mathcal{B}=\widehat{S}\cup V_{12}\cup D_1\cup D_2 \cup B^3$, where $B^3$ denotes a $3$-ball, by the irreducibility of $K(\alpha)$.
Let $H_1(\widehat{S}\cup V_{12})=\langle l,m,x\rangle=\mathbf{Z}\oplus \mathbf{Z}\oplus\mathbf{Z}$, where
$l,m$ are the cycles on $\widehat{S}$ as in Figure \ref{worst1} and $x$ is represented by the core of $V_{12}$ directed from $u_1$ to $u_2$.
Then
$$H_1(\mathcal{B})=\langle l,m,x\rangle / \langle \partial D_1,\partial D_2 \rangle,$$
and $\partial D_1=2x+l, \partial D_2=2m$ with suitable orientations.

Similarly, consider the bigon $E_1$ between $d_1$ and $b_1$ and the $6$-gon face $E_2$ bounded by $c_1,c_3,g,a_4,a_2,e$ in $G_T$.
Then we have $\mathcal{W}=\widehat{S}\cup V_{21}\cup E_1\cup E_2\cup B^3$.
Thus 
$$H_1(\mathcal{W})=\langle l,m,y\rangle / \langle \partial E_1,\partial E_2 \rangle,$$
where $y$ is represented by the core of $V_{21}$ directed from $u_2$ to $u_1$,
and $\partial E_1=2y-l-m, \partial E_2=2l+4m$ with suitable orientations.

Hence $H_1(K(\alpha))=\langle l,m,x,y\rangle / \langle \partial D_1,\partial D_2,\partial E_1, \partial E_2 \rangle=\mathbf{Z}_4\oplus \mathbf{Z}_4$,
which is not cyclic.
This is a contradiction.
\end{proof}

Thus we have shown that the case $p_1=t/2$ is impossible.


\subsection*{Case 2. $p_1=t/2+1$}

Since $p_2+p_3+p_4+p_5=4t-2$, 
at least two of $p_i$ are $t$.
By the parity rule, $p_2+p_3$ and $p_4+p_5$ are even.
Thus we may assume that $(p_2,p_3,p_4,p_5)=(t,t,t,t-2)$ or $(t,t,t-1,t-1)$ without loss of generality.
Let $A, B, C$ and $D$ be the families of parallel negative edges of $G_S$ with $p_2, p_3, p_4$ and $p_5$ edges, respectively.

\begin{lemma}
$(p_2,p_3,p_4,p_5)=(t,t,t-1,t-1)$ is impossible.
\end{lemma}

\begin{proof}
Assume $(p_2,p_3,p_4,p_5)=(t,t,t-1,t-1)$.
Each edge of $A$ has labels with the same parity at its ends by the parity rule.
Then any edge of $C$ has labels with opposite parities at its ends.
This contradicts the parity rule.
\end{proof}

\begin{figure}[tb]
\includegraphics*[scale=0.4]{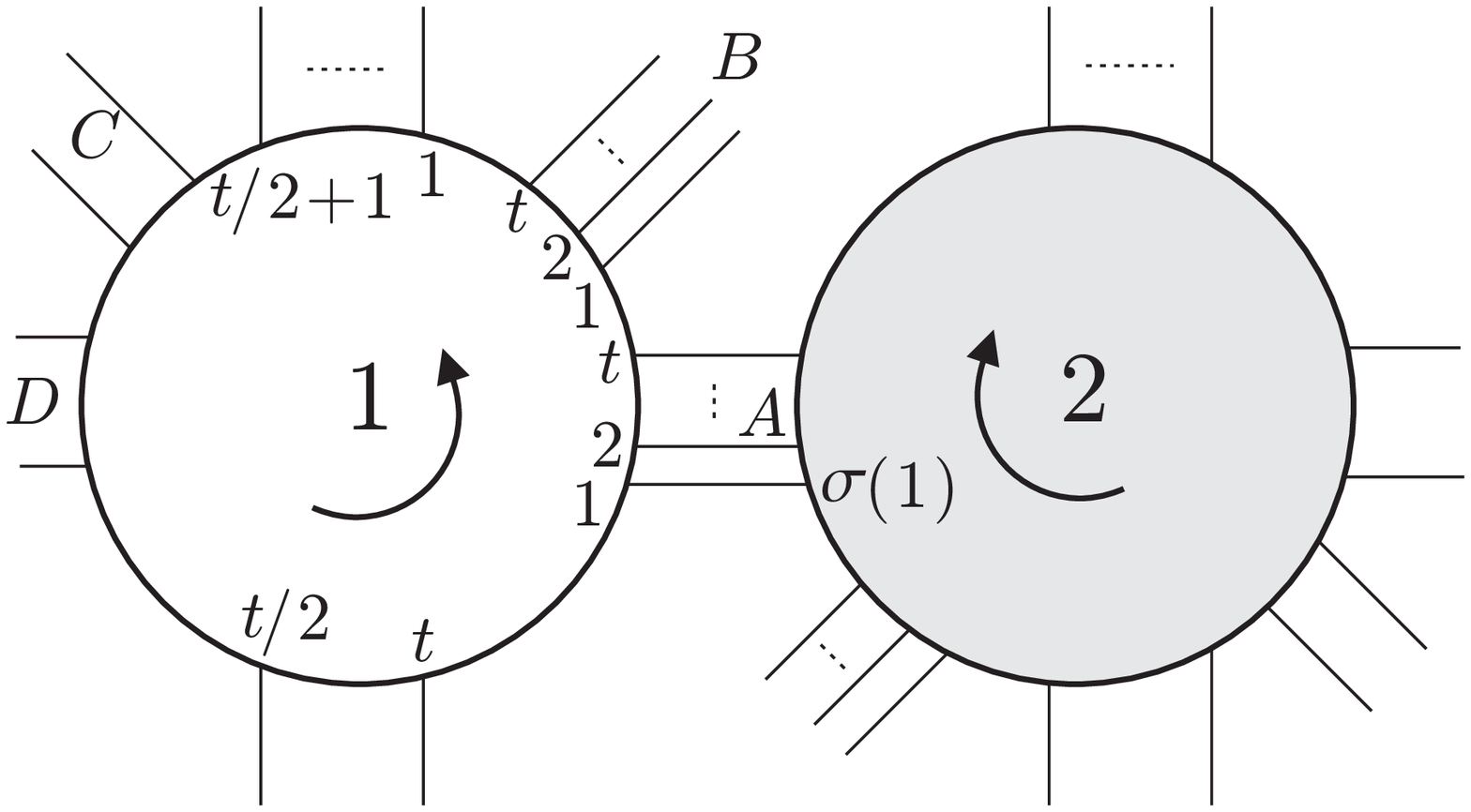}
\caption{}\label{s2case2}
\end{figure}

Thus we have $(p_2,p_3,p_4,p_5)=(t,t,t,t-2)$.
The labels in $G_S$ can be assumed as shown in Figure \ref{s2case2}.
Let $\sigma$ be the permutation associated to $A$ as before.
Then there is an $S$-cycle with label pair $\{t/2,t/2+1\}$ among positive loops at vertex $u_1$.

\begin{lemma}\label{ABD}
$\sigma$ is the identity.
\end{lemma}

\begin{proof}
Among the positive loops at vertex $u_2$, there is an $S$-cycle with label pair $\{\sigma(1)+t/2-1,\sigma(1)+t/2\}$.
(All labels are considered modulo $t$.)
Assume $\sigma(1)\ne 1$.
Then $\sigma(1)\ge 3$ and odd.

If $t=4$, then $\sigma(1)=3$, and hence $G_S$ has two $S$-cycles with label pairs $\{2,3\}$ and $\{4,1\}$.
This is impossible by Lemma \ref{Scharlemann}(2).

Assume $t>4$.
Then we see that $\sigma^2$ is the identity and $\sigma(1)=t/2+1$ by Lemma \ref{involution}.
(The argument applies here without change.)
Hence $G_S$ has two $S$-cycles with label pairs $\{t/2,t/2+1\}$ and $\{t,1\}$ respectively. 
The edges of $A$ form cycles of length two on $\widehat{T}$, and there are at least four such cycles.
In particular, $v_1$ and $v_{t/2+1}$ lie on the same cycle, and so do $v_{t/2}$ and $v_{t}$.
But we cannot locate the edges of the above two $S$-cycles to satisfy Lemma \ref{Schsupp} simultaneously.
\end{proof}

\begin{lemma}
$t=4$.
\end{lemma}

\begin{proof}
Assume $t>4$.
Then we see that the edges of $C$ form two essential cycles on $\widehat{T}$.
By Lemma \ref{ABD}, each vertex of $G_T$ is incident to a loop.
Thus there would be a trivial loop.
\end{proof}

\begin{lemma}
$t=4$ is impossible.
\end{lemma}

\begin{proof}
In $G_T$, $v_1$ and $v_4$ are incident to $3$ loops, and $v_2$ and $v_3$ are incident to two loops.
In $G_S$, there are two $S$-cycles with label pair $\{2,3\}$.
The edges of them give four edges between $v_2$ and $v_3$ in $G_T$.
Then two endpoints with label $1$ of loops at $v_2$ cannot be successive among the five occurrences of label $1$,
which contradicts Lemma \ref{jumping}.
\end{proof}

Hence the case $p_1=t/2+1$ is also impossible.

\section{The case that $s=t=2$}\label{s2t2}

Finally, we consider the case that $s=t=2$.
Then both $\overline{G}_S$ and $\overline{G}_T$ are subgraphs of the graph in Figure \ref{t2}.
If $K(\gamma)$ contains a Klein bottle, then $\gamma$ is a multiple of four \cite{Te}.
Hence either $K(\alpha)$ or $K(\beta)$ does not contain a Klein bottle, because $\Delta(\alpha,\beta)=|\alpha-\beta|=5$.
Without loss of generality, we can assume that $K(\beta)$ does not contain a Klein bottle.
Also, we use the notation $p_i$ for the number of
edges in the families of parallel positive or negative edges in $G_S$ as in the previous section.

\begin{lemma}\label{max}
$1\le p_1\le 3$ and $p_i\le 2$ for $i=2,3,4,5$.
\end{lemma}

\begin{proof}
If $p_1\ge 4$, then there are two bigons among loops which lie on the same side of $\widehat{T}$.
By Lemma \ref{bothparallel}, the four edges of the bigons belong to mutually distinct families of parallel negative edges in $G_T$.
But this implies that $K(\beta)$ contains a Klein bottle (\cite[The proof of Lemma 5.2]{G-L2}), a contradiction. 

If a family of parallel negative edges contains three edges in $G_S$, then
two of them are incident to the same vertex in $G_T$.
Thus they are also parallel in $G_T$, which contradicts Lemma \ref{bothparallel}.
Hence $p_i\le 2$ for $i\ne 1$.

Since $u_1$ has degree $10$ in $G_S$, $2p_1+p_2+p_3+p_4+p_5=10$.
Thus $p_1\ge 1$.
\end{proof}

\begin{lemma}\label{1loop}
$p_1=1$ is impossible.
\end{lemma}

\begin{proof}
If $p_1=1$, then $p_i=2$ for $i\ne 1$ by Lemma \ref{max}.
Then $G_T$ has the same form as in Figure \ref{s2case1id}.
But the jumping number argument in the proof of Lemma \ref{notid} eliminates these configurations again.
\end{proof}

\begin{lemma}\label{2loop}
$p_1=2$ is impossible.
\end{lemma}

\begin{figure}[tb]
\includegraphics*[scale=0.65]{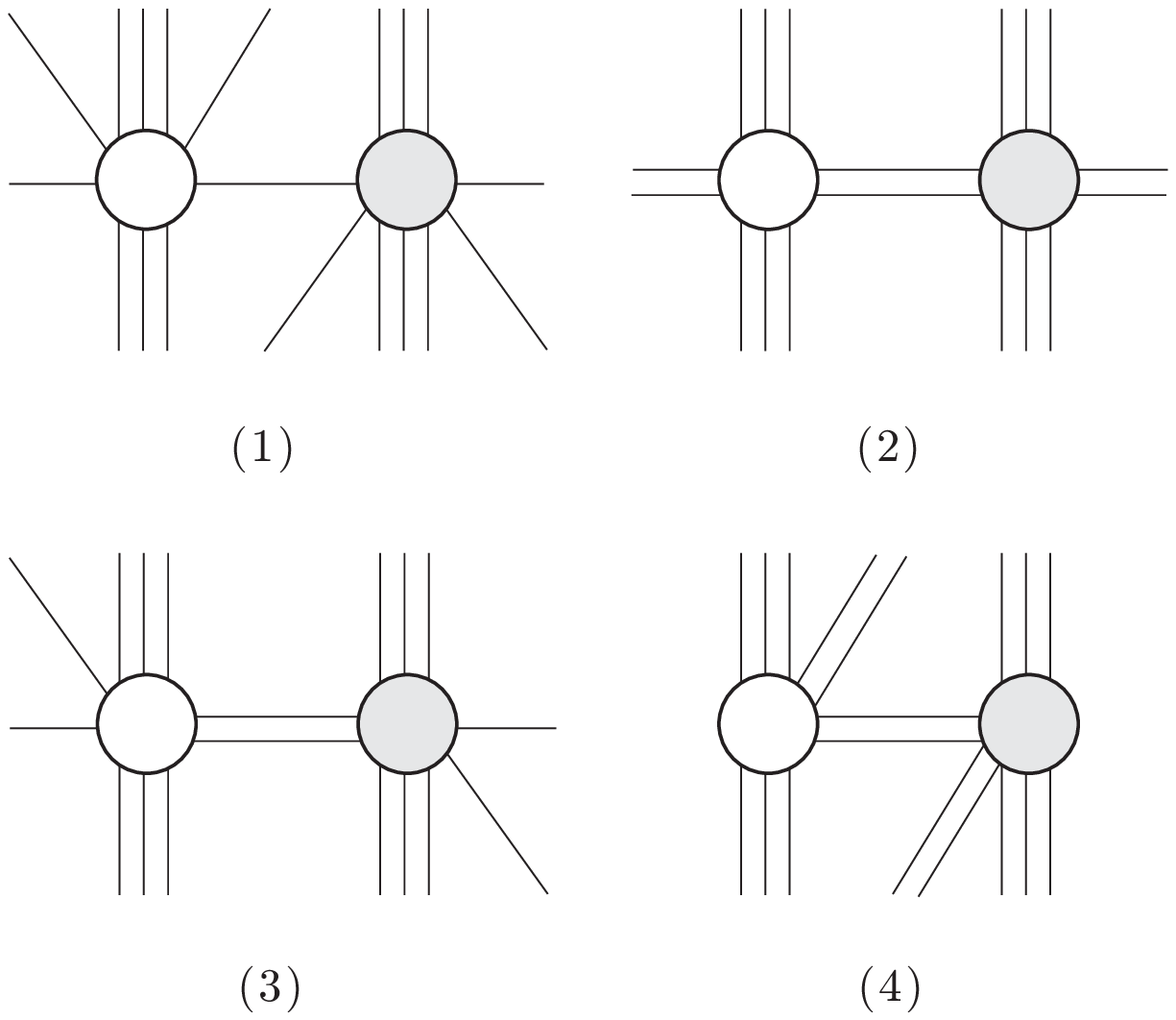}
\caption{}\label{s2t2case2}
\end{figure}

\begin{figure}[tb]
\includegraphics*[scale=0.8]{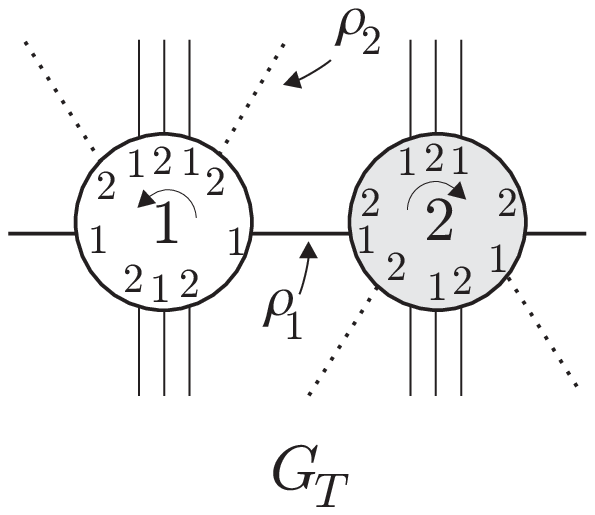}
\caption{}\label{s2t2case2irr}
\end{figure}

\begin{figure}[tb]
\includegraphics*[scale=0.7]{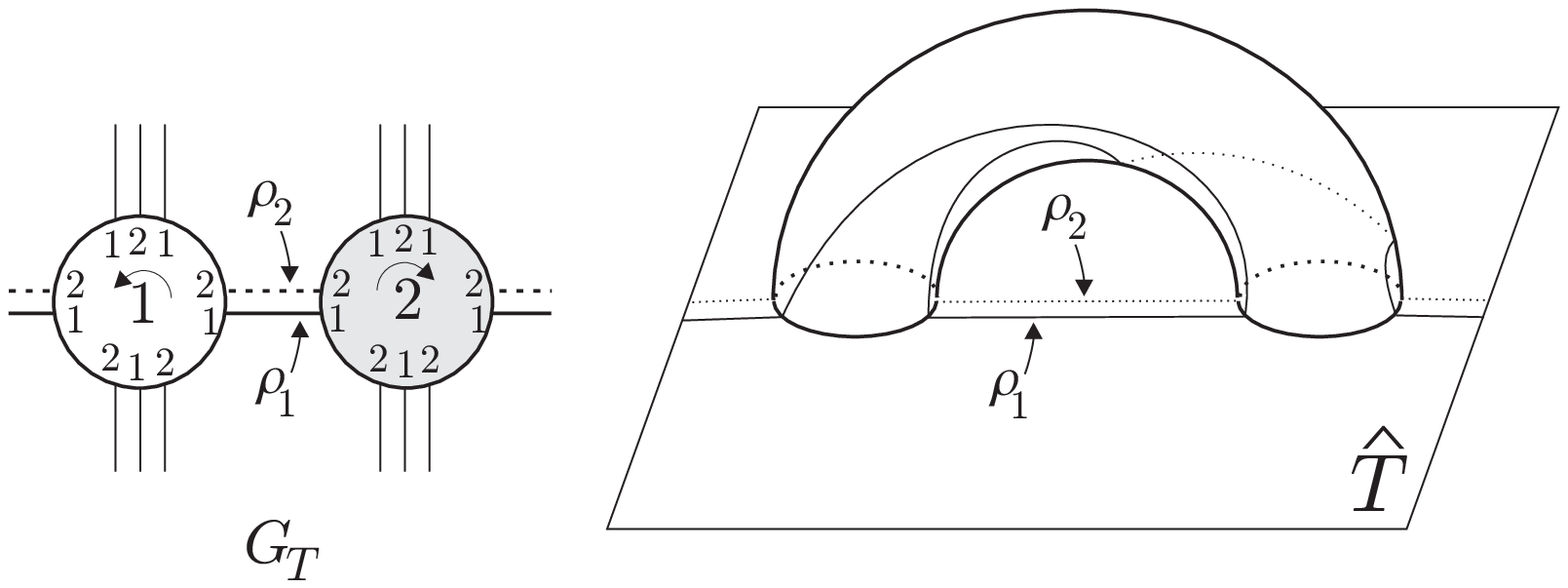}
\caption{}\label{s2t2case2irr0}
\end{figure}

\begin{proof}
Assume $p_1=2$.  Then we can assume that $p_2+p_3=2$ and $p_4+p_5=4$ by the parity rule.
Hence $p_4=p_5=2$.
If $p_2=p_3=1$, then the labels in $G_S$ contradicts the parity rule.
Thus we can assume that $p_2=2$ and $p_3=0$.
Then there are $4$ possibilities for $G_T$ as in Figure \ref{s2t2case2}.
(As in the proof of Lemma \ref{max}, a family of parallel negative edges in $G_T$ contains
at most two edges.) 
We see that (3) contradicts the parity rule.

(4) can be eliminated by a jumping number argument.
In $G_T$, there are two negative edges incident to $v_1$ with the same label $1$.
Their endpoints at $v_1$ are consecutive among the five occurrences of label $1$.
But these points are not consecutive at $u_1$, which contradicts Lemma \ref{jumping}.

To eliminate (1) and (2), note that
$G_S$ contains two $S$-cycles $\rho_1$ and $\rho_2$ whose
faces lie on the same side of $\widehat{T}$.
From the labeling of $G_T$, we can determine the edges of $\rho_i$ in $G_T$
as in Figure \ref{s2t2case2irr} for (1) and Figure \ref{s2t2case2irr0} for (2).
In the former, $K(\beta)$ contains a Klein bottle as in the proof of Lemma \ref{max}, a contradiction.
In the latter, it is impossible to connect these edges on $\partial V_\beta$ simultaneously (see Figure \ref{s2t2case2irr0}).
\end{proof}

\begin{lemma}\label{3loop}
$p_1=3$ is impossible.
\end{lemma}

\begin{proof}
Assume $p_1=3$.
Then we can assume that $(p_2+p_3,p_4+p_5)=(4,0)$ or $(2,2)$.
If $(p_2+p_3,p_4+p_5)=(4,0)$ then $p_2=p_3=2$.
The endpoints of two negative edges with label $1$ are successive at $u_1$ among the five occurrences of label $1$.
By Lemma \ref{jumping}, those points are also successive at $v_1$ among the five occurrences of label $1$.
Then if we put six negative edges between $v_1$ and $v_2$, then
there would be a pair of edges which is parallel in both graphs, a contradiction by Lemma \ref{bothparallel}.

If $(p_2+p_3,p_4+p_5)=(2,2)$, then there are three possibilities for $G_S$ as in Figure \ref{s2t2case2}(1), (2) and (3).
Then (3) contradicts the parity rule.
If $G_S$ is (2),
then the labeling of $G_S$ implies that $G_T$ has two parallel loops at each vertex.
Thus $G_T$ has two $S$-cycles.
It is easy to see that their faces lie on the same side of $\widehat{S}$.
Hence the argument in the proof of Lemma \ref{2loop} works again (with an exchange of roles between $G_S$ and $G_T$).

For (1), $G_S$ and $G_T$ are determined as shown in Figure \ref{s2t2final}.
Then we can conclude that $K(\alpha)$ contains a Klein bottle as in the proof of Lemma \ref{2loop},
but this is not a contradiction.

\begin{figure}[tb]
\includegraphics*[scale=0.7]{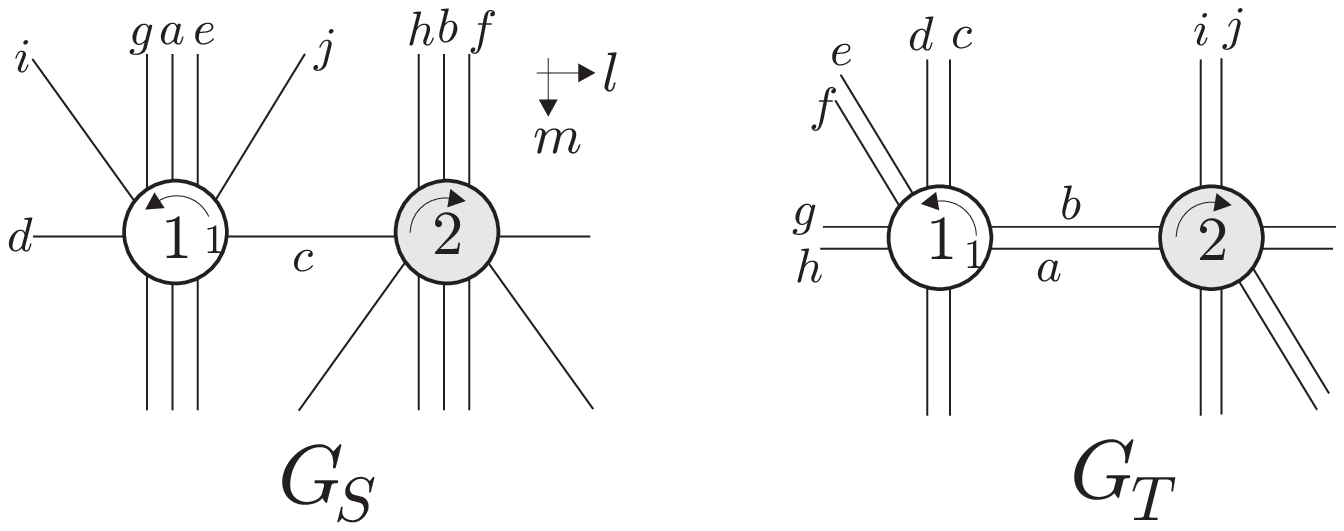}
\caption{}\label{s2t2final}
\end{figure}

To eliminate (1), we calculate $H_1(K(\alpha))$ and $H_1(K(\beta))$.
First, the (unique) edge correspondence between $G_S$ and $G_T$ is shown in Figure \ref{s2t2final}.
Let $D_1$ and $D_2$ be the bigons between the edges $c$ and $d$, $a$ and $b$, respectively, in $G_T$.
Also, let $E_1$ and $E_2$ be the $3$-gon bounded by $d,e,h$ and the $4$-gon bounded by $a,b,c,i$ in $G_T$, respectively. 
As in the proof of Lemma \ref{lem:worst2}, let us call $\mathcal{B}$ the side of $\widehat{S}$ which contains $D_1$ and $D_2$,
and call $\mathcal{W}$ the other side.
Then we can see that
$\mathcal{B}=\widehat{S}\cup V_{12}\cup D_1\cup D_2\cup B^3$ and $\mathcal{W}=\widehat{S}\cup V_{21}\cup E_1\cup E_2\cup B^3$,
where $V_{12}=V_\alpha\cap \mathcal{B}$ and $V_{21}=V_\alpha\cap\mathcal{W}$ as in the proof of Lemma \ref{lem:worst2}.
Hence 
\begin{eqnarray*}
H_1(\mathcal{B}) &=& \langle l,m,x \rangle / \langle \partial D_1,\partial D_2 \rangle = \langle l,m,x \rangle / \langle 2x+l,2m \rangle, \cr
H_1(\mathcal{W}) &=& \langle l,m,y \rangle / \langle \partial E_1,\partial E_2 \rangle = \langle l,m,y \rangle / \langle y-2m-l,3m+l \rangle,
\end{eqnarray*}
where $H_1(\widehat{S})=\langle l,m \rangle$ (see Figure \ref{s2t2final}), and $x$ and $y$ are represented by the cores of $V_{12}$ and $V_{21}$
directed from $u_1$ to $u_2$, and from $u_2$ to $u_1$, respectively.
By the Mayer-Vietoris sequence, 
$$H_1(K(\alpha))=\langle l,m,x,y \rangle / \langle \partial D_1,\partial D_2,\partial E_1,\partial E_2 \rangle =\mathbf{Z}_4.$$
This means that $|\alpha|=4$.

Similarly, we calculate $H_1(K(\beta))$.
Let $D_1'$ be the bigon between $a$ and $e$, and $D_2'$ the $3$-gon bounded by $c,h,j$ in $G_S$. 
Also, let $E_1'$ be the bigon between $a$ and $g$, and $E_2'$ be the $3$-gon bounded by $c,e,j$.
By using these, we can build up $K(\beta)$ as
$$K(\beta)=\widehat{T}\cup V_{\beta} \cup D_1'\cup D_2'\cup E_1'\cup E_2'\cup (\mbox{two}\ 3\mbox{-balls}).$$
Then we can show that $H_1(K(\beta))=\mathbf{Z}_{11}$, which implies $|\beta|=11$.
This contradicts that $|\alpha-\beta|=5$.
\end{proof}

By Lemmas \ref{1loop}, \ref{2loop} and \ref{3loop}, the case $s=t=2$ is impossible.
Hence the proof of Theorem \ref{main1} is now complete.

\section{Questions}

We ask some questions:
\begin{itemize}
\item[(1)] If a hyperbolic knot admits three toroidal slopes, then is it either the figure-eight knot or the $(-2,3,7)$-pretzel knot?
\item[(2)] If a hyperbolic knot has two integral toroidal slopes with distance $4$, then does at least one toroidal surgery yield a Klein bottle?
\end{itemize}


\section*{Acknowledgements}
The author would like to thank the referee for valuable suggestions and comments, and
Mario Eudave-Mu\~{n}oz for pointing out an error in the original manuscript.

\bibliographystyle{amsplain}

\end{document}